\documentclass[12pt]{article}
\usepackage{latexsym,amssymb,amsmath}
\begin{document}
\baselineskip=18pt

\newcommand{\la}{\langle}
\newcommand{\ra}{\rangle}
\newcommand{\psp}{\vspace{0.4cm}}
\newcommand{\pse}{\vspace{0.2cm}}
\newcommand{\ptl}{\partial}
\newcommand{\dlt}{\delta}
\newcommand{\sgm}{\sigma}
\newcommand{\al}{\alpha}
\newcommand{\be}{\beta}
\newcommand{\G}{\Gamma}
\newcommand{\gm}{\gamma}
\newcommand{\vs}{\varsigma}
\newcommand{\Lmd}{\Lambda}
\newcommand{\lmd}{\lambda}
\newcommand{\td}{\tilde}
\newcommand{\vf}{\varphi}
\newcommand{\for}{\mbox{for}}
\newcommand{\wt}{\mbox{wt}\:}
\newcommand{\der}{\mbox{Der}}
\newcommand{\ad}{\mbox{ad}\:}
\newcommand{\stl}{\stackrel}
\newcommand{\ol}{\overline}
\newcommand{\ul}{\underline}
\newcommand{\es}{\epsilon}
\newcommand{\dmd}{\diamond}
\newcommand{\clt}{\clubsuit}
\newcommand{\mbb}{\mathbb}
\newcommand{\llra}{\Longleftrightarrow}
\newcommand{\vt}{\vartheta}
\newcommand{\rta}{\rightarrow}
\newcommand{\ves}{\varepsilon}
\newcommand{\dg}{\dag}

To appear in Journal of Lie Theory {\bf 16} (2006), no. 4.
\vspace{1cm}

\begin{center}{\Large \bf Tree Diagram Lie Algebras of  Differential Operators}\end{center}
\begin{center}{\Large \bf and Evolution Partial Differential Equations}\end{center} \vspace{0.2cm}

\begin{center}{\large Xiaoping Xu}\footnote{Research Supported
 by China NSF 10431040}
\end{center}

\begin{center}{Institute of mathematics, Academy of Mathematics \&
System Sciences}\end{center}
\begin{center}{Chinese Academy of Sciences, Beijing 100080, P. R.
China}\end{center}

\vspace{2cm}

 \begin{center}{\Large\bf Abstract}\end{center}
\vspace{1cm} {\small A tree diagram is a tree with positive
integral weight on each edge, which is a notion generalized from
the Dynkin diagrams of finite-dimensional simple Lie algebras. We
introduce two nilpotent Lie algebras and their extended solvable
Lie algebras associated with each tree diagram. The solvable tree
diagram Lie algebras turn out to be complete Lie algebras of
maximal rank analogous to the Borel subalgebras of
finite-dimensional simple Lie algebras. Their abelian ideals are
completely determined. Using a high-order Campbell-Hausdorff
 formula and certain abelian ideals of the tree diagram Lie algebras, we solve the initial value
 problem of first-order evolution partial differential equations
 associated with nilpotent tree diagram Lie
algebras and high-order evolution partial differential equations,
including heat conduction type equations related to generalized
Tricomi operators associated with trees.} \vspace{0.8cm}

\section{Introduction}

In 1905, Schur [Si] proved that the maximum number of linearly
independent commuting $n\times n$ matrices is $[|n^2/4|]+1$.
Jacobson [J1] (1944) gave a simplified derivation of Schur's
result. Malcev [M] (1945) determined the abelian subalgebras of
maximal dimension of a finite-dimensional semisimple Lie algebra.
Kostant [K1] (1965) found a connection of Malcev's result with the
maximal eigenvalue of the Laplacian acting on the exterior algebra
of the Lie algebra via the adjoint representation. Moreover, he
[K2] (1998) related abelian ideals of a Borel subalgebra to the
discrete series representations of the Lie group. Panyushev and
R\"{o}hrle [PR] (2001) described spherical orbits in terms of
abelian ideals. Furthermore, Panyushev [P] (2003) gave a
correspondence of maximal abelian ideals of a Borel subalgebra to
long positive roots. Suter [Sr] (2004) determined the maximal
dimension among abelian subalgebras of a finite-dimensional simple
Lie algebra purely in terms of certain invariants and gave a
uniform explanation for Panyushev's result. Besides, he gave a
generalization and explanation of the symmetry property of a
certain subposet of Young's lattice. Kostant [K3] (2004) showed
that the powers of the Euler product and abelian ideals of a Borel
subalgebra are intimately related. There were many results on
non-abelian ideals of Borel subalgebras ( cf. [AKOP], [CP1-CP4],
[KOP]) and structure of nilpotent Lie algebras (e.g., cf. [FS],
[GK], [Sl1], [Sl2]).

Barros-Neto and  Gel'fand [BG1,BG2] (1998,2002) studied
fundamental solutions for the Tricomi operator. Our initial
motivation of this work is to solve heat conduction type evolution
partial differential equation related to generalized Tricomi
operators associated with trees. It turns out that some abelian
ideals of certain ``tree diagram solvable Lie algebras" have to be
used. This motivates us to introduce two nilpotent Lie algebras
and their extended solvable Lie algebras associated with a general
tree diagram. Familiar examples of such solvable Lie algebras are
the Borel subalgebras of special linear Lie algebras and
symplectic Lie algebras. The abelian ideals of these solvable Lie
algebras are completely determined. Moreover, we completely solve
the first-order evolution partial differential equations
 associated with nilpotent tree diagram Lie
algebras and high-order evolution partial differential equations,
including heat conduction type equations related to generalized
Tricomi operators associated with trees. Below we give a more
detailed technical introduction.

Barros-Neto and Gel'fand [BG1,BG2] (1998, 2002) studied  solutions
of the equation
$$u_{xx}+xu_{yy}=\dlt(x-x_0,y-y_0)\eqno(1.1)$$
related to the Tricomi operator $\ptl_x^2+x\ptl_y^2$. A natural
generalization of the Tricomi operator is
$\ptl_{x_1}^2+x_1\ptl_{x_2}^2+\cdots+x_{n-1}\ptl_{x_n}^2$. The
equation
$$u_t=\ptl_{x_1}^2(u)+\ptl_{x_2}^2(u)+\cdots+\ptl_{x_n}^2(u)\eqno(1.2)$$
is a well known classical heat conduction equation related to the
Laplacian operator
$\ptl_{x_1}^2+\ptl_{x_2}^2+\cdots+\ptl_{x_n}^2$. As pointed out in
[BG1, BG2], the Tricomi operator is an analogue of the Laplacian
operator. An immediate  analogue of heat conduction equation is
$$u_t=\ptl_{x_1}^2(u)+x_1\ptl_{x_2}^2(u)+\cdots+x_{n-1}\ptl_{x_n}^2(u).\eqno(1.3)$$
Denote by $\mbb{C}$ the field of complex numbers. Graphically, the
above equation is related to the Dynkin diagram of the special
linear Lie algebra $sl(n+1,\mbb{C})$:

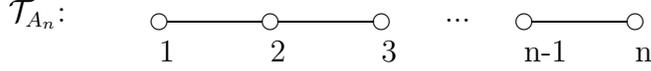
\setlength{\unitlength}{3pt}

\begin{picture}(70,8)\put(2,0){${\cal T}_{A_n}$:}
\put(21,0){\circle{2}}\put(22,0){\line(1,0){12}}\put(21,-5){1}
\put(35,0){\circle{2}}\put(35,-5){2} \put(36,0){\line(1,0){12}}
\put(49,0){\circle{2}}\put(49,-5){3}
\put(57,0){...}\put(67,0){\circle{2}} \put(67,-5){n-1}
\put(68,0){\line(1,0){12}}\put(81,0){\circle{2}}\put(81,-5){n}
\end{picture}
\vspace{0.8cm}

\noindent Naturally, we should also consider similar equations
related to the graph:

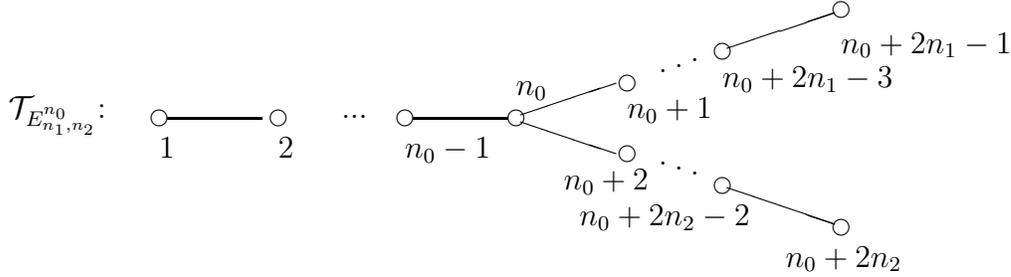
\begin{picture}(108,25)\put(2,0){${\cal T}_{E^{n_0}_{n_1,n_2}}$:}
\put(21,0){\circle{2}}\put(21,-5){1}\put(22,0){\line(1,0){12}}
\put(36,0){\circle{2}}\put(36,-5){2}
\put(44,0){...}\put(52,0){\circle{2}}\put(52,-5){$n_0-1$}\put(53,0){\line(1,0){12}}
\put(66,0){\circle{2}}\put(66,3){$n_0$}\put(66.6,0.5)
{\line(3,1){12}}\put(80,4.5){\circle{2}}\put(80,0){$n_0+1$}\put(84,5.8){.}
\put(86,6.46){.}\put(88,7.14){.}\put(92,8.5){\circle{2}}\put(92,4){$n_0+2n_1-3$}
\put(92.6,8.8){\line(3,1){13.6}}
\put(107,13.7){\circle{2}}\put(107,8.3){$n_0+2n_1-1$}
\put(66.6,-0.5)
{\line(3,-1){12}}\put(80,-4.5){\circle{2}}\put(72,-9){$n_0+2$}\put(84,-5.8){.}
\put(86,-6.46){.}\put(88,-7.14){.}\put(92,-8.5){\circle{2}}\put(74,-13.5){$n_0+2n_2-2$}
\put(92.6,-8.8){\line(3,-1){13.6}}
\put(107,-13.7){\circle{2}}\put(100,-18.7){$n_0+2n_2$}
\end{picture}
\vspace{2.5cm}

\noindent which is the Dynkin diagram of the orthogonal Lie
algebra $so(2n,\mbb{C})$ when $n_1=n_2=1$, and the Dynkin diagram
of the simple Lie algebra of types $E_6,\;E_7,\;E_8$ if
$(n_0,n_1,n_2)=(3,1,2),\;(3,1,3),\;(3,1,4)$, respectively. In
fact, when $(n_0,n_1,n_2)=(3,2,2),\;(4,1,3),\;(6,1,2)$, it is also
the Dynkin diagram of the affine Kac-Moody Lie algebra of types
$E_6^{(1)},\;E_7^{(1)},\;E_8^{(1)}$,  respectively (cf. [Kv]).
These diagrams are special examples of trees in graph theory.

A {\it tree} ${\cal T}$ consists of a finite set of {\it nodes}
${\cal N}= \{\iota_1,\iota_2,...,\iota_n\}$ and a set of {\it
edges}
$${\cal E}\subset\{(\iota_i,\iota_j)\mid 1\leq i<j\leq n\}\eqno(1.4)$$
such that for each node $\iota_i\in{\cal N}$, there exists a
unique sequence $\{\iota_{i_1},\iota_{i_2},...,\iota_{i_r}\}$ of
nodes with $1=i_1<i_2<\cdots<i_{r-1}<i_r=i$ for which
$$(\iota_{i_1},\iota_{i_2}),(\iota_{i_2},\iota_{i_3}),...,(\iota_{i_{r-2}},\iota_{i_{r-1}}),
(\iota_{i_{r-1}},\iota_{i_r})\in{\cal E}.\eqno(1.5)$$  We also
denote the tree ${\cal T}=({\cal N},{\cal E})$. We  identify a
tree ${\cal T}=({\cal N},{\cal E})$ with a graph by depicting a
small circle for each node in ${\cal N}$ and a segment connecting
$i$th circle to $j$th circle for the edge
$(\iota_i,\iota_j)\in{\cal E}$ (cf. the above Dynkin diagrams of
type $A$ and $E$).

 For a tree ${\cal T}=({\cal N},{\cal E})$, we call the
differential operator
$$d_{\cal T}=\ptl_{x_1}^2+\sum_{(\iota_i,\iota_j)\in{\cal
E}}x_i\ptl_{x_j}^2\eqno(1.6)$$ a {\it generalized Tricomi operator
of type} ${\cal T}$. Moreover, we call the partial differential
equation
$$u_t=d_{\cal T}(u)\eqno(1.7)$$
a {\it generalized heat conduction  equation associated with the
tree} ${\cal T}$, where $u$ is a function in $t,x_1,x_2,...,x_n$.
For instance, the generalized heat equation of type ${\cal
T}_{E^{n_0}_{n_1,n_2}}$ is:

\begin{eqnarray*}\hspace{1.5cm}u_t&=&(\ptl_{x_1}^2+\sum_{i=1}^{n_0-1}x_i\ptl_{x_{i+1}}^2
+\sum_{r=0}^{n_2-1}x_{n_0+2r}\ptl_{x_{n_0+2r+2}}^2\\ &&
+x_{n_0}\ptl_{x_{n_0+1}}^2 +\sum_{i=1}^{n_1-1}x_{n_0+2i-1}
\ptl_{x_{n_0+2i+1}}^2)(u).\hspace{5.3cm}(1.8)\end{eqnarray*}

Let $m_1,m_2,...,m_n$ be $n$ positive integers. The difficulty of
solving the equation (1.7) is the same as that of solving the
following more general evolution partial differential equation:
$$u_t=(\ptl_{x_1}^{m_1}+\sum_{(\iota_i,\iota_j)\in{\cal
E}}x_i\ptl_{x_j}^{m_j})(u).\eqno(1.9)$$ One of our main objectives
is to solve the equation (1.9) subject to the initial condition:
$$u(0,x_1,...,x_n)=f(x_1,...,x_n)\qquad\for\;\;x_i\in[-a_i,a_i],\eqno(1.10)$$
where $f$ is a given continuous function and $a_i$ are given
positive real constants.

Our strategy of solving the above problem is first to attack
simpler  first-order evolution partial differential equations.
Denote by $\mbb{Z}_+$ the set of positive integers. We define a
tree diagram ${\cal T}^d$ to be a tree ${\cal T}=({\cal N},{\cal
E})$ with a weight map $d:{\cal E}\rta \mbb{Z}_+$, denoted as
${\cal T}^d=({\cal N},{\cal E},d)$. Moreover, we associate a tree
diagram ${\cal T}^d$ with a first-order partial differential
equation:
$$u_t=(\ptl_{x_1}+\sum_{(\iota_i,\iota_j)\in{\cal
E}}x_i^{d[(\iota_i,\iota_j)]}\ptl_{x_j})(u).\eqno(1.11)$$ The
above equation is a ``dual" of (1.9) in a certain sense. It is
also an equation of polynomial transformations for vector fields
and may have practical applications. Our another objective is to
solve (1.11) subject to the initial value condition:
$$u(0,x_1,x_2,....,x_n)=f(x_1,x_2,...,x_n)\eqno(1.12)$$
for any given first-order differentiable function $f$.

We remark that completely solving a linear partial differential
equation with large number of variables is in general a difficult
problem. The method of characteristic lines usually works only for
linear partial differential equations with few variables. In order
to solve (1.11), we have to use certain hidden symmetry.

Given a positive integer $n$, we denote
$$\mbb{A}=\sum_{m_1,m_2,...,m_n=0}^{\infty}\mbb{C}[x_1,x_2,...,x_n]
\ptl_{x_1}^{m_1}\ptl_{x_2}^{m_2}\cdots
\ptl_{x_n}^{m_n},\eqno(1.13)$$ the associative algebra of
differential operators in $n$ variables. Moreover, we define a Lie
bracket on $\mbb{A}$ by
$$[A,B]=AB-BA\qquad\for\;\;A,B\in \mbb{A}.\eqno(1.14)$$
For any tree diagram ${\cal T}^d=({\cal N},{\cal E},d)$, we define
the {\it upward nilpotent Lie algebra of differential operators}
associated with ${\cal T}^d$ by
$$L_0({\cal T}^d)=\mbox{the Lie subalgebra of}\;\mbb{A}\;\mbox{generated
by}\:\{\ptl_{x_1},x_i^{d[(\iota_i,\iota_j)]}\ptl_{x_j}\mid
(\iota_i,\iota_j)\in{\cal E}\}.\eqno(1.15)$$ A node $\iota_i$ of a
tree ${\cal T}$ is called a {\it tip} if there does not exist
$i\leq j\leq n$ such that $(\iota_i,\iota_j)\in{\cal E}$. Set
$$\Psi=\mbox{the set of all tips of}\;{\cal T}.\eqno(1.16)$$
The {\it downward nilpotent Lie algebra of differential operators}
associated with ${\cal T}^d$ is defined by
\begin{eqnarray*}\hspace{2cm}{\cal L}_0({\cal T}^d)&=&\mbox{the
Lie subalgebra of}\;\mbb{A}\;\mbox{generated by}\\ &
&\{\ptl_{x_r},x_j^{d[(\iota_i,\iota_j)]}\ptl_{x_i}\mid
\iota_r\in\Psi,\;(\iota_i,\iota_j)\in{\cal
E}\}.\hspace{3.8cm}(1.17)\end{eqnarray*} The Lie algebra
$L_0({\cal T}^d)$ is the hidden symmetry of the equation (1.11),
while the hidden symmetry of the equation (1.9) is of the form
${\cal L}_0(\ol{\cal T}^d)$ for some other tree diagram $\ol{\cal
T}^d$. Both $L_0({\cal T}^d)$ and ${\cal L}_0({\cal T}^d)$ are
nilpotent Lie algebras of maximal rank (e.g., cf. [Sl1], [GK]).

Denote $$H=\sum_{i=1}^n\mbb{C}x_i\ptl_{x_i}.\eqno(1.18)$$ Extend
the above nilpotent Lie algebras to:
$$L_1({\cal
T}^d)=H+L_0({\cal T}^d),\qquad {\cal L}_1({\cal T}^d)=H+{\cal
L}_0({\cal T}^d),\eqno(1.19)$$ respectively.  According to [MZJ],
both $L_1({\cal T}^d)$ and ${\cal L}_1({\cal T}^d)$ are complete
solvable Lie algebras. Indeed, $L_1({\cal T}_{A_n}^d)$ is
isomorphic to a Borel subalgebra of $sl(n+1,\mbb{C})$ if $d({\cal
E})=\{1\}$ and isomorphic to a Borel subalgebra of
$sp(2n,\mbb{C})$ when $d({\cal
E}\setminus\{(\iota_{n-1},\iota_n)\})=\{1\}$ and
$d[(\iota_{n-1},\iota_n)]=2$. The Lie algebras $L_1({\cal T}^d)$
and ${\cal L}_1({\cal T}^d)$  are natural analogues of the Borel
subalgebras of finite-dimensional simple Lie algebras.

In Section 2, we will completely determine the abelian ideals of
the Lie algebra $L_1({\cal T}^d)$. The abelian ideals of the Lie
algebra ${\cal L}_1({\cal T}^d)$ will be determined in Section 3.
In Section 4, we will use some abelian ideals of $L_1({\cal T}^d)$
and a high-order Campbell-Hausdorff formula to solve the equation
(1.11) subject to the condition (1.12). In Section 5, we will use
the results in Section 4 to solve the equation (1.9) subject to
the condition (1.10).

We would like to thank the referee for helpful comments.

\section{Upward Tree Diagram Lie Algebras}

In this section, we will study the structure of the Lie algebra
$L_0({\cal T}^d)$ in (1.16), mainly the abelian ideals of
$L_1({\cal T}^d)$ in (1.19).

For a node $\iota_i$ in a tree ${\cal T}$, the unique sequence
$${\cal C}_i=\{\iota_{i_1},\iota_{i_2},...,\iota_{i_r}\}\eqno(2.1)$$ of nodes with
$1=i_1<i_2<\cdots<i_{r-1}<i_r=i$ satisfying (1.5) is called the
{\it clan} of the node $\iota_i$. The node $\iota_{i_{r-1}}$ is
called the {\it parent} of the node $\iota_i$, denoted as
$\iota_{p(i)}$. It can be verified that the following set
\begin{eqnarray*}B({\cal
T}^d)&=&\{\ptl_{x_1},x_{i_1}^{s_2d[(\iota_{i_1},\iota_{i_2})]-s_1}
x_{i_2}^{s_3d[(\iota_{i_2},\iota_{i_3})]-s_2} \cdots
x_{i_{r-1}}^{d[(\iota_{i_{r-1}},\iota_{i_r})]-s_{r-1}}\ptl_{x_i}\\
& &\mid 2\leq i\leq n,\; {\cal
C}_i=\{\iota_{i_1},\iota_{i_2},...,\iota_{i_r}\},\;0\leq
s_{r-1}\leq d[(\iota_{i_{r-1}},\iota_{i_r})],\\ & &0\leq s_\es\leq
s_{\es+1}d[(\iota_{i_\es},\iota_{i_{\es+1}})]\}.\hspace{8cm}(2.2)\end{eqnarray*}
is a basis of $L_0({\cal T}^d)$.

In order to find the dimension formula of $L_0({\cal T}^d)$, we
define $\ell(m_1)=m_1+1$ and \begin{eqnarray*} \hspace{1cm}& &
\ell(m_1,m_2,...,m_r,m)=\mbox{the coefficient of}\;t^{m_1\cdots
m_rm}\;\\ & &\mbox{in}\;
\frac{1}{(1-t)^2(1-t^{m_1})(1-t^{m_1m_2})\cdots (1-t^{m_1m_2\cdots
m_r})}\hspace{4.5cm}(2.3)\end{eqnarray*} for any positive integers
$m_1,\cdots,m_r,m$. For instance,
$$\ell(1,1,...,1,m)={r+1+m\choose{m}}.\eqno(2.4)$$
Denote by $\mbb{N}$ the additive semigroup of nonnegative
integers. Then
\begin{eqnarray*}\hspace{1cm}& &\{x_1^{i_2m_1-i_1}x_2^{i_3m_2-i_2}\cdots
x_{r-1}^{i_rm_{r-1}-i_{r-1}}x_r^{im_r-i_r}x_{r+1}^{m-i}\\ & &\mid
0\leq i\leq m,\;0\leq i_r\leq im_r,\;0\leq i_\es\leq
i_{\es+1}m_\es\}\\ &=&\{x_1^{s_0}x_2^{s_1}\cdots
x_r^{s_{r-1}}x_{r+1}^{s_r}\mid
(s_0,s_1,...,s_r)\in\mbb{N}^{\:r+1},\\ &
&s_0+\sum_{i=1}^rs_im_1m_2\cdots m_i\leq m_1m_2\cdots
m_rm\}.\hspace{5.7cm}(2.5)\end{eqnarray*} The cardinality of (2.5)
is exactly $\ell(m_1,...,m_r,m)$. \psp

{\bf Proposition 2.1}. {\it We have:
$$\mbox{\it dim}\:L_0({\cal
T}^d)=1+\sum_{i=2}^n\ell(d[(\iota_{i_1},\iota_{i_2})],...,
d[(\iota_{i_{r-1}},\iota_{i_r})]),\eqno(2.6)$$ where we have used
the notion (2.1).}\psp

{\bf Example 2.1}. If ${\cal T}={\cal T}_{A_n}$ (diagram after
(1.3)) with $d({\cal E}\setminus \{(\iota_{n-1},\iota_n)\})=\{1\}$
and $d[(\iota_{n-1},\iota_n)]=m$, then
$$\dim L_0({\cal T}^d_{A_n})={n+m-1\choose
m}+\frac{n(n-1)}{2}.\eqno(2.7)$$

When ${\cal T}={\cal T}_{E^{n_0}_{n_1,n_2}}$ (diagram before
(1.4)) with
$$d({\cal E}\setminus \{(\iota_{n_0+2n_1-3},\iota_{n_0+2n_1-1}),
(\iota_{n_0+2n_2-2},\iota_{n_0+2n_2})\})=\{1\}\eqno(2.8)$$ and
$$d[(\iota_{n_0+2n_1-3},\iota_{n_0+2n_1-1})]=m_1,\;\;
d[(\iota_{n_0+2n_2-2},\iota_{n_0+2n_2})]=m_2,\eqno(2.9)$$ then
\begin{eqnarray*}\dim L_0({\cal T}^d_{A_n})&=&{n_0+n_1+m_1-1\choose m_1}
+{n_0+n_2+m_2-1\choose m_2}\\
&&+\frac{n_0^2+n_1^2+n_2^2-3n_0-n_1-n_2}{2}+n_0(n_1+n_2).\hspace{3.7cm}(2.10)\end{eqnarray*}
\vspace{0.1cm}

 For two subspaces $U$ and $V$ of a Lie algebra ${\cal G}$, we
define
$$[U,V]=\mbox{Span}\:\{[u,v]\mid u\in U,\;v\in V\}.\eqno(2.11)$$
The central series
$${\cal G}^1={\cal G},\qquad{\cal G}^i=[{\cal G},{\cal
G}^{i-1}].\eqno(2.12)$$ The Lie algebra ${\cal G}$ is called {\it
nilpotent} if there exists a positive integer $m$ such that ${\cal
G}^m=\{0\}$. For a nilpotent Lie algebra ${\cal G}$, the
nilpotence $N({\cal G})$ is the smallest positive integer $m$ such
that ${\cal G}^m=\{0\}$. In fact, our Lie algebra $L_0({\cal
T}^d)$ is nilpotent. Suppose that the clan of a node $\iota_i$ is
(2.1). We define
$$N(i)=1+\sum_{s=2}^rd[(\iota_{i_1},\iota_{i_2})]d[(\iota_{i_2},\iota_{i_3})]\cdots
d[(\iota_{i_{s-1}},\iota_{i_s})].\eqno(2.13)$$ Recall the notion
$\Psi$ in (1.16).\psp

{\bf Proposition 2.2}. {\it We get}:
$$N(L_0({\cal T}^d))=\mbox{\it max}\;\{N(i)\mid \iota_i\in\Psi\}.\eqno(2.14)$$ \vspace{0.1cm}

The Lie algebra $L_0({\cal T}^d)$ is in fact a nilpotent Lie
algebra of maximal rank. \psp

{\bf Proposition 2.3}. {\it We have}:
$$\mbox{\it the center of}\;L_0({\cal
T}^d)=\sum_{\iota_i\in\Psi}\mbb{C}\ptl_{x_i}.\eqno(2.15)$$
\vspace{0.1cm}

Note that $H=\sum_{i=1}^n\mbb{C}x_i\ptl_{x_i}$ is  a Cartan
subalgebra of the Lie algebra $L_1({\cal T}^d)=H+L_0({\cal T}^d)$.
Denote by $H^\ast$ the space of linear functions on $H$. Define
$\ves_i\in H^\ast$ by
$$\ves_i(x_j\ptl_{x_j})=\dlt_{i,j}.\eqno(2.16)$$
A nonzero vector $u\in\mbb{A}$ is a {\it root vector} if there
exists $0\neq \al\in H^\ast$ such that
$$[h,u]=\al(h)u\qquad\for\;h\in H.\eqno(2.17)$$
The function $\al$ is called a {\it root}. Observe that
$$x_1^{m_1}\cdots x_i^{m_i}\ptl_{x_{i+1}}\;\mbox{is a root vector
corresponding to the
root}\;\sum_{r=1}^im_r\ves_r-\ves_{i+1}.\eqno(2.18)$$ Expressions
(2.2) and (2.16)-(2.18) imply:\psp

{\bf Lemma 2.4}. {\it The subspace $L_0({\cal T}^d)$ forms an
ideal of $L_1({\cal T}^d)$ and is a direct sum of root subspaces
of dimension 1 with respect to $H$. Any abelian ideal of
$L_1({\cal T}^d)$ is a direct sum of some root subspaces of
$L_0({\cal T}^d)$}.\psp

By (2.2), the set of roots of $L_0({\cal T}^d)$:
\begin{eqnarray*} \qquad& &R(L_0({\cal
T}^d))=\{-\ves_1,\sum_{s=1}^{r-1}j_s\ves_s-\ves_i\mid 2\leq i\leq
n,\;{\cal C}_i=\{\iota_{i_1},\iota_{i_2},...,\iota_{i_r}\},\\&
&j_1+\sum_{s=2}^{r-1}j_s d[(\iota_{i_1},\iota_{i_2})]\cdots
d[(\iota_{i_{s-1}},\iota_{i_s})]\leq
d[(\iota_{i_1},\iota_{i_2})]\cdots
d[(\iota_{i_{r-1}},\iota_{i_r})]\}\hspace{1.6cm}(2.19)\end{eqnarray*}
with respect to $H$.

For convenience, we denote
$$v_{-\ves_1}=\ptl_{x_1},\;\;v_{_{\sum_{s=1}^{r-1}j_s\ves_s-\ves_i}}=x_{i_1}^{j_1}\cdots
x_{i_{r-1}}^{j_{r-1}}\ptl_{x_i}.\eqno(2.20)$$ Suppose that $I$ is
a nonzero abelian ideal of $L_1({\cal T}^d)$. Define the set of
roots of $I$ by:
$$R(I)=\{\al\in R(L_0({\cal
T}^d))\mid v_\al\in I\}.\eqno(2.21)$$ Then
$$I=\mbox{Span}\:\{v_\al\mid \al\in R(I)\}.\eqno(2.22)$$ Determining $I$ is equivalent to determining
all possible $R(I)$.\psp

{\bf Lemma 2.5}. {\it Let $2\leq i\leq n$ and ${\cal
C}_i=\{\iota_{i_1},\iota_{i_2},...,\iota_{i_r}\}$. If
$d[(\iota_{i_{r-1}},\iota_i)]>1$, then}
$$-\ves_1,j_1\ves_{i_1}+\cdots+ j_{s-1}\ves_{i_{s-1}}-\ves_s\not\in R(I)\qquad \for\;\;2\leq s<r.\eqno(2.23)$$

{\it Proof}. Suppose that (2.23) fails. If $-\ves_1\in R(I)$, then
$\ptl_{x_1}\in I$. So
$$[\ptl_{x_1},x_1\ptl_{x_{i_2}}]=\ptl_{x_{i_2}}\in I.\eqno(2.24)$$ Assume $j_1\ves_{i_1}+\cdots
j_{s-1}\ves_{i_{s-1}}-\ves_s\in R(I)$ for some $2\leq s<r$,
equivalently,
$$x_{i_1}^{j_1}\cdots x_{i_{s-1}}^{j_{s-1}}\ptl_{x_{i_s}}\in I.\eqno(2.25)$$ Hence
$$(\mbox{ad}\:\ptl_{x_{i_1}})^{j_1}\cdots(\mbox{ad}\:\ptl_{x_{i_{s-1}}})^{j_{s-1}}
(x_{i_1}^{j_1}\cdots
x_{i_{s-1}}^{j_{s-1}}\ptl_{x_{i_s}})=j_1!\cdots
j_{s-1}!\ptl_{x_{i_s}}\in I.\eqno(2.26)$$ Thus we  can always
assume $\ptl_{x_{i_s}}\in I$ for some $2\leq s<r$.

Note
$$\ptl_{x_{i_{r-1}}}=[...[\ptl_{x_{i_s}},x_{i_s}\ptl_{x_{i_{s+1}}}],...,x_{i_{r-2}}
\ptl_{x_{i_{r-1}}}]\in I.\eqno(2.27)$$ On the other hand,
$x_{i_{r-1}}^2\ptl_{x_i}\in L_0({\cal T}^d)$. Thus
$$x_{i_{r-1}}\ptl_{x_i}=\frac{1}{2}[\ptl_{x_{i_{r-1}}},x_{i_{r-1}}^2\ptl_{x_i}]\in I.\eqno(2.28)$$ But
$$0\neq
\ptl_{x_i}=[\ptl_{x_{i_{r-1}}},x_{i_{r-1}}\ptl_{x_i}],\eqno(2.29)$$
which contradicts the fact that $I$ is abelian. $\qquad\Box$\psp

For any $2\leq i\leq n$ and ${\cal
C}_i=\{\iota_{i_1},\iota_{i_2},...,\iota_{i_r}\}$, we define
\begin{eqnarray*}\qquad R_i&=&\{(j_1,j_2,...,j_{r-1})\in\mbb{N}^{\:r-1}\mid
j_1+\sum_{s=2}^{r-1}j_s d[(\iota_{i_1},\iota_{i_2})]\cdots
d[(\iota_{i_{s-1}},\iota_{i_s})]\\ & &\leq
d[(\iota_{i_1},\iota_{i_2})]\cdots
d[(\iota_{i_{r-1}},\iota_{i_r})]\}.\hspace{7.3cm}(2.30)\end{eqnarray*}
Denote
$$x^{\vec j}=x_{i_1}^{j_1}x_{i_2}^{j_2}\cdots
x_{i_{r-1}}^{j_{r-1}}\qquad\for\;\;\vec j=(j_1,j_2,...,j_{r-1})\in
R_i.\eqno(2.31)$$ Moreover, we define a partial ordering
``$\prec$'' on $R_i$ by
$$(j_1,j_2,...,j_{r-1})\prec(l_1,l_2,...,l_{r-1})\eqno(2.32)$$
if $$j_s+\sum_{\es=s+1}^{r-1}j_\es
d[(\iota_{i_s},\iota_{i_{s+1}})]\cdots
d[(\iota_{i_{\es-1}},\iota_{i_{\es}})]\leq
l_s+\sum_{\es=s+1}^{r-1}l_\es
d[(\iota_{i_s},\iota_{i_{s+1}})]\cdots
d[(\iota_{i_{\es-1}},\iota_{i_{\es}})]\eqno(2.33)$$ for any $1\leq
s\leq r-1$. For convenience, we denote $R_1=\{0\}$ and $x^0=1$.

A node $\iota_i$ of a tree diagram ${\cal T}^d$ is called
multiplicity-free if there exists $j<i$ such that
$(\iota_j,\iota_i)\in{\cal E}$ and $d[(\iota_j,\iota_i)]=1$. A
node $\iota_j$ is called a {\it descendant} of $\iota_i$ if $i<j$
and $\iota_i\in{\cal C}_j$. Set
$${\cal D}_i=\mbox{the set of all descendants
of}\;\iota_i.\eqno(2.34)$$ Denote
$$\Upsilon=\{\iota_i\in{\cal N}\mid \iota_j\;\mbox{is
multiplicity-free}\;\mbox{for any}\;\iota_j\in{\cal
D}_i\}.\eqno(2.35)$$
 A subset $S$ of nodes is called an independent if each element of
$S$ is not a descendent of any other element in $S$. For any
independent subset $S$ of $\Upsilon$, we define
$$I(S)=\mbox{Span}\:\{x^{\vec j}\ptl_{x_i},x^{\vec
j}\ptl_{x_s}\mid\iota_i\in S,\;\vec j\in R_i,\;\iota_s\in{\cal
D}_i\}.\eqno(2.36)$$ It can be verified that $I(S)$ is an abelian
ideal of $L_1({\cal T}^d)$. \psp

{\bf Theorem 2.6}. {\it Any maximal abelian ideal of $L_1({\cal
T}^d)$ is of the form $I(S)$ for some maximal independent subset
$S$ of $\Upsilon$.}

{\it Proof.} Suppose that $I$ is an abelian ideal of $L_1({\cal
T}^d)$. Denote
$$W=\{\iota_i\in{\cal N}\mid x^{\vec j}\ptl_{x_i}\in
I\;\mbox{for some}\;\vec j\in R_i\}.\eqno(2.37)$$ By Lemma 2.5,
$W\subset\Upsilon$. Define
$$S=\{\iota_i\in W\mid \iota_i\not\in {\cal D}_r\;\for\;\iota_r\in
W\setminus\{\iota_i\}\}.\eqno(2.38)$$ Then $S$ is an independent
subset of $\Upsilon$.

Any $\iota_s\in W\setminus S$ must be a descendent of some
$\iota_i\in S$, that is, there exist
$$i_1=i<i_2<\cdots<i_{r-1}<i_r=s\;\;\mbox{with}\;\;(\iota_\es,\iota_{\es+1})\in{\cal
E}.\eqno(2.39)$$ Since $\iota_i\in W$, we have
$$\ptl_{x_{i_\es}}\in I\qquad\for\;\;1\leq\es\leq
r\eqno(2.40)$$ by the proof of Lemma 2.5. Suppose $x^{\vec
j}\ptl_{x_s}\in I$ for some $\vec j\in R_s$. Then
$$0=[\ptl_{x_{i_\es}},x^{\vec
j}\ptl_{x_s}]=\ptl_{x_{i_\es}}(x^{\vec j})\ptl_{x_s}\qquad
\for\;\;1\leq\es\leq r.\eqno(2.41)$$ By Lemma 2.5,
$$d[(\iota_{i_\es},\iota_{i_{\es+1}})]=1\qquad\for\;\;1\leq\es\leq
r-1.\eqno(2.42)$$ Expressions (2.31), (2.41) and (2.42) imply
$$x^{\vec j}=x^{\vec l}\qquad\mbox{for some}\;\;\vec l\in
R_i.\eqno(2.43)$$ Hence, $x^{\vec j}\ptl_{x_s}\in I(S)$. That is,
$I\subset I(S)$.

Suppose that $I$ is a maximal abelian ideal. Then $I=I(S)$ and $S$
is a maximal independent subset $S$ of $\Upsilon.\qquad\Box$ \psp

{\bf Corollary 2.7}. {\it If all the tips are not
multiplicity-free, then $I(\Psi)$ is the unique maximal abelian
ideal of $L_1({\cal T}^d)$}.\psp

Let $I$ be any nonzero abelian ideal of $L_1({\cal T}^d)$. Define
$W$ and $S$ as in (2.37) and (2.38). We have already proved
 $I\subset I(S)$. Recall the ordering defined in
(2.32) and (2.33). Given $\iota_i\in{\cal N}$ and $\vec j\in R_i$,
we define
$$I[i,\vec j]=\mbox{Span}\:\{x^{\vec l}\ptl_{x_i},x^{\vec
l}\ptl_{x_r}\mid r\in{\cal D}_i,\;\vec j\succ\vec l\in
R_i\}.\eqno(2.44)$$ If $i\in\Upsilon$, $I[i,\vec j]$ is an abelian
ideal of $L_1({\cal T}^d)$ generated by $x^{\vec j}\ptl_{x_i}$.

Given $i\in S$, we define
$$\Pi_i=\{\vec j\in R_i\mid x^{\vec j}\ptl_{x_i}\in I\},
\qquad \Pi_r=\{\vec l\in R_i\mid x^{\vec l}\ptl_{x_r}\in{\cal
I}\}\eqno(2.45)$$ for $r\in{\cal D}_i$. Moreover, we set
$$K_i=\mbox{the set of all maximal elements
in}\;\Pi_i.\eqno(2.46)$$ For any $\iota_r\in{\cal D}_i$, the
parent node  $\iota_{p(r)}$ is the unique node  such that
$(\iota_{p(r)},\iota_r)\in{\cal E}$. So $p(r)<r$, and $p(r)=i$ or
$\iota_{p(r)}\in {\cal D}_i$. In fact, $\Pi_{p(r)}\subset \Pi_r$.
Now we define
$$K_r=\mbox{the set of all maximal elements
in}\;\Pi_r\setminus\Pi_{p(r)}\qquad\for\;\;r\in{\cal
D}_i.\eqno(2.47)$$  Since
$$W=S\bigcup(\bigcup_{\iota_i\in S}{\cal D}_i),\eqno(2.48)$$
the set $$\{x^{\vec j}\ptl_{x_r}\mid r\in W,\;\vec j\in
K_r\}\eqno(2.49)$$ is the unique (up to constant multiples)
minimal set of root vectors that generate $I$. We call the pair
$(S,\{K_r\mid \iota_r\in W\})$ an {\it admissible pair for abelian
ideal.} Note that given any $r\in W$ and $\iota_s\in {\cal D}_r$,
$$\vec l\not\prec \vec j\qquad \for\;\;\vec j\in K_r,\;\vec l\in
K_s.\eqno(2.50)$$

Conversely, we can construct admissible pairs as follows. A subset
$K$ of $R_i$ is called {\it independent} if $K$ itself is the set
of maximal elements in $K$. Pick a non-empty independent set
$S\subset \Upsilon$. Define $W$ by (2.48). For each $\iota_i\in S$
and $\iota_\es\in{\cal D}_i$, we choose independent subsets
$K_i,K_\es\subset R_i$ such that $K_i\neq\emptyset$ and (2.50)
holds. We define
$$I[S,\{K_r\mid \iota_r\in W\}]=\sum_{\iota_s\in
W}\sum_{\vec j\in K_s}I[s,\vec j]\eqno(2.51)$$ (cf. (2.44)). Then
$I[S,\{K_r\mid \iota_r\in W\}]$ is a nonzero abelian ideal with
$(S,\{K_r\mid \iota_r\in W\})$ as the corresponding admissible
pair. This proves:\psp

 {\bf Theorem 2.8}. {\it Any nonzero abelian ideal of
$L_1({\cal T}^d)$ must be of the form $I[S,\{K_r\mid \iota_r\in
W\}]$ for some admissible pair $(S,\{K_r\mid \iota_r\in W\})$
constructed in the above paragraph. So there is a one-to-one
correspondence between the set of nonzero abelian ideals of
$L_1({\cal T}^d)$ and the set of admissible pairs.}\psp

For any two integers $i$ and $j$ with $i\leq j$, we use the
following notion of indices in the rest of this paper:
$$\ol{i,j}=\{i,i+1,...,j-1,j\}.\eqno(2.52)$$
\vspace{0.1cm}

{\bf Example 2.2}. Let ${\cal T}={\cal T}_{A_n}$  and $d({\cal
E})=\{1\}$. The Lie algebra $L_1({\cal T}_{A_n}^d)$ is isomorphic
to a Borel subalgebra of $sl(n+1,\mbb{C})$. In this case,
$\Upsilon={\cal N}$. The number of maximal abelian ideals of
$L_1({\cal T}_{A_n}^d)$ is $n$. Denote
$$\hat{0}=(0,...,0),\;\;\hat{i}=(0,...,0,\stl{i}{1},0,...,0).\eqno(2.53)$$
Then
$$R_i=\{\hat{0},\hat{1},...,\widehat{i-1}\}\eqno(2.54)$$
and
$$\hat{0}\prec\hat{1}\prec\cdots\prec\widehat{i-1}.\eqno(2.55)$$
Any nonempty independent set $S$ of $\Upsilon$ consists of one
element, say, $S=\{\iota_i\}$. Now ${\cal
D}_i=\{\iota_{i+1},\iota_{i+2},...,\iota_n\}$. So
$$W=\{\iota_i,\iota_{i+1},\iota_{i+2},...,\iota_n\}.\eqno(2.56)$$
Any nonempty $K_r$ must also consist of one element in $R_i$.
Write
$$\{\iota_r\in W\mid
K_r\neq\emptyset\}=\{\iota_{i_1}=\iota_i,\iota_{i_2},...,\iota_{i_s}\}\eqno(2.57)$$
with $i_1<i_2<\cdots<i_s$. Then
$$K_{i_r}=\{\hat{j_r}\},\qquad j_r\in\ol{0,i-1}.\eqno(2.58)$$
Condition (2.53) is equivalent to
$$j_1<j_2<\cdots <j_{s-1}<j_s.\eqno(2.59)$$
So there exists a one-to-one correspondence between the set of
admissible pairs and the set of sequences:
$$j_1<j_2<\cdots <j_{s-1}<j_s<i<i_2<\cdots i_s\eqno(2.60)$$
in $\ol{0,n}$. But the sequence (2.60) is equivalent to a subset
of the set $\ol{0,n}$.

Let the empty set corresponding to the zero ideal of $L_1({\cal
T}^d_{A_n})$. Thus the number of abelian ideals $L_1({\cal
T}^d_{A_n})$ is equal to the number of subsets with even number of
elements in  $\ol{0,n}$. So there are exactly $2^n$ abelian ideals
in $L_1({\cal T}^d_{A_n})$. This result was obtained by Peterson,
 using the corresponding affine Weyl
group (mentioned in [K3], [Sr]). \psp

{\bf Example 2.3}. Consider ${\cal T}={\cal T}_{A_n}$ with
$d({\cal E}\setminus\{(\iota_{n-1},\iota_n)\})=\{1\}$ and
$d[(\iota_{n-1},\iota_n)]=2$. The Lie algebra $L_1({\cal
T}_{A_n}^d)$ is isomorphic to a Borel subalgebra of
$sp(2n,\mbb{C})$. In this case, $\Upsilon=\{\iota_n\}$. So
$L_1({\cal T}_{A_n}^d)$ has a unique maximal abelian ideal. Treat
$$\hat{r}=\hat{0}+\hat{r}.\eqno(2.61)$$ We have
$$R_n=\{\hat{i}+\hat{j}\mid 0\leq i\leq j\leq
n-1\}.\eqno(2.62)$$

Note for $i_1\leq i_2$ and $j_1\leq j_2$,
$$\hat{i}_1+\hat{i}_2\prec \hat{j}_1+\hat{j}_2\llra  i_1\leq
j_1\;\mbox{and}\;i_2\leq j_2.\eqno(2.63)$$ A nonempty independent
set of $R_n$ is either of the form $\{\hat{i}_s+\hat{j}_s\mid
s\in\ol{1,r}\}$ with
$$0\leq i_1<i_2<\cdots < i_r<j_r<j_{r-1}<\cdots <j_1\leq
n-1\eqno(2.64)$$ or of the form
$\{\hat{i}_s+\hat{j}_s,2\hat{j}_{r+1}\mid s\in\ol{1,r}\}$ with
$$0\leq i_1<i_2<\cdots < i_r<j_{r+1}<j_r<j_{r-1}<\cdots <j_1\leq
n-1.\eqno(2.65)$$ Thus there exists a one-to-one correspondence
between the set of nonzero abelian ideals and the set of nonempty
subsets of $\ol{0,n-1}$. Therefore, the number of abelian ideals
of $L_1({\cal T}^d_{A_n})$ is $2^n$, which was proved by Peterson,
using the corresponding affine Weyl group (mentioned in [K3],
[Sr]).\psp

{\bf Example 2.4}. Suppose ${\cal T}={\cal T}_{A_3}$,
$d[(\iota_1,\iota_2)]=1$ and $d[(\iota_2,\iota_3)]=m$, where $m>2$
is an integer. Now $\Upsilon=\{\iota_3\}$ and
$$R_3=\{(i,j)\in\mbb{N}^{\:2}\mid i+j\leq m\}.\eqno(2.66)$$
Moreover,
$$(i_1,j_1)\prec (i_2,j_2)\llra i_1+j_1\leq
i_2+j_2\;\mbox{and}\;j_1\leq j_2.\eqno(2.67)$$ So any nonempty
independent subset of $R_3$ must be of the form
$\{(i_1,j_1),(i_2,j_2),...,(i_r,j_r)\}$ with
$$0\leq j_1<j_2<\cdots<j_r\leq
i_r+j_r<i_{r-1}+j_{r-1}<\cdots< i_1+j_1\leq m.\eqno(2.68)$$ This
shows the number of abelian ideals in $L_1({\cal T}^d_{A_3})$ is
exactly $2^{m+1}$.\psp

{\bf Example 2.5}. Suppose that all nodes of ${\cal T}^d$ except
tips are multiplicity-free and
$$d[(\iota_{p(r)},\iota_r)]=2\qquad\for\;\;\iota_r\in\Psi\eqno(2.69)$$
(cf. (1.16)). By Corollary 2.7, $L_1({\cal T}^d)$ has a unique
maximal abelian ideal. Moreover, Example 2.3 implies
$$\mbox{the number of abelian ideals of}\;L_1({\cal
T}^d)=2^{\sum_{\iota_i\in\Psi}|{\cal C}_i|}.\eqno(2.70)$$
\vspace{0.1cm}

{\bf Example 2.6}. Let ${\cal T}={\cal T}_{E^{n_0}_{n_1,n_2}}$
with $n_0>1$ and $n_1,n_2\geq 1$. Assume
$$d({\cal
E}\setminus\{(\iota_{n_0+2n_1-3},\iota_{n_0+2n_1-1})\}=\{1\},\;\;d[(\iota_{n_0+2n_1-3},
\iota_{n_0+2n_1-1})]=2.\eqno(2.71)$$ Then
$$\Upsilon=\{\iota_{n_0+2n_1-1},\iota_{n_0+2i}\mid
i\in\ol{1,n_2}\}.\eqno(2.72)$$

 The number of maximal abelian ideal of $L_1({\cal
T}_{E^{n_0}_{n_1,n_2}}^d)$ is $n_2$. By Examples 2.2 and 2.3, the
number of abelian ideals is
$$2^{n_0+n_1}\left(2^{n_2}+\sum_{i=1}^{n_2}\sum_{r=1}^i{i\choose r}{n_0\choose r}
2^{n_2-i}\right).\eqno(2.73)$$

\section{Downward Tree Diagram Lie Algebras}

In this section, we will study the structure of the Lie algebra
${\cal L}_0({\cal T}^d)$ in (1.17), mainly the abelian ideals of
${\cal L}_1({\cal T}^d)$ in (1.19).

Recall the downward nilpotent Lie algebra of differential
operators associated with a tree diagram ${\cal T}^d$ is
\begin{eqnarray*}\hspace{2cm}{\cal L}_0({\cal T}^d)&=&\mbox{the Lie
subalgebra of}\;\mbb{A}\;\mbox{generated by}\\ &
&\{\ptl_{x_r},x_j^{d[(\iota_i,\iota_j)]}\ptl_{x_i}\mid
(\iota_i,\iota_j)\in{\cal
E},\;\iota_r\in\psi\}\hspace{4.2cm}(3.1)\end{eqnarray*} (cf.
(1.16)). To find a basis for ${\cal L}_0({\cal T}^d)$, we set
$${\cal F}_i=\{\iota_i\}\bigcup {\cal D}_i,\;\;{\cal
C}_{i,j}=\{\iota_i\}\bigcup ({\cal C}_j\setminus {\cal
C}_i)\eqno(3.2)$$ for $\iota_i\in{\cal N}$ and $\iota_j\in{\cal
D}_i$. Moreover, we denote
$${\cal E}_i=\{(\iota_r,\iota_s)\in{\cal E}\mid
\iota_r,\iota_s\in{\cal F}_i\},\;\;{\cal
E}_{i,j}=\{(\iota_r,\iota_s)\in{\cal E}\mid
\iota_r,\iota_s\in{\cal C}_{i,j}\},\eqno(3.3)$$  and let
$$\kappa_i=\prod_{(\iota_r,\iota_s)\in{\cal
E}_i}d[(\iota_r,\iota_s)],\qquad
\kappa_{i,j}=\frac{\kappa_i}{\prod_{(\iota_r,\iota_s)\in{\cal
E}_{i,j}}d[(\iota_r,\iota_s)]}.\eqno(3.4)$$ By (2.5),
 the following set
$${\cal B}({\cal
T}^d)=\{\ptl_{x_\es},(\prod_{\iota_s\in{\cal D}_i} x^{j_s}_s)
\ptl_{x_i}\mid
\iota_\es\in\Psi,\;j_s\in\mbb{N},\;\sum_{\iota_s\in{\cal
D}_i}j_s\kappa_{i,s}\leq \kappa_i\}\eqno(3.5)$$ is a basis of
${\cal L}_0({\cal T}^d)$.

In order to find the dimension formula of ${\cal L}_0({\cal
T}^d)$, we define
$$\be_i=\mbox{the coefficient
of}\;t^{\kappa_i}\;\mbox{in}\;\frac{1}{(1-t)\prod_{\iota_s\in{\cal
D}_i}(1-t^{\kappa_{i,s}})}.\eqno(3.6)$$ \vspace{0.1cm}

{\bf Proposition 3.1}. {\it We have:
$$\mbox{\it dim}\:{\cal L}_0({\cal
T}_d)=|\Psi|+\sum_{\iota_i\in{\cal
N}\setminus\Psi}\be_i.\eqno(3.7)$$ Moreover, the nilpotence}
$$N({\cal L}_0({\cal T}^d))=1+\sum_{\iota_j\in{\cal D}_1}\;\prod_{(\iota_r,\iota_s)\in{\cal
E}_{1,j}}d[(\iota_r,\iota_s)].\eqno(3.8)$$ \vspace{0.1cm}

The Lie algebra ${\cal L}_0({\cal T}^d)$ is in fact a nilpotent
Lie algebra of maximal rank. \psp

{\bf Proposition 3.2}. {\it We have}:
$$\mbox{\it the center of}\;{\cal L}_0({\cal
T}^d)=\mbb{C}\ptl_{x_1}.\eqno(3.9)$$ \vspace{0.1cm}

Propositions 2.2 and 3.2 imply that ${\cal L}_0({\cal T}^d)$ can
not be isomorphic to $L_0(\ol{\cal T}^d)$ for any other tree
diagram $\ol{\cal T}^d$ if ${\cal T}$ has more than one tips.

Note that $H=\sum_{i=1}^n\mbb{C}x_i\ptl_{x_i}$ is a toral Cartan
subalgebra of the Lie algebra ${\cal L}_1({\cal T}^d)=H+{\cal
L}_0({\cal T}^d)$  and ${\cal L}_0({\cal T}^d)$ is the set of all
ad-nilpotent elements of ${\cal L}_1({\cal T}^d)$.

According to (2.16)-(2.18), the operator
$$x_{i+1}^{m_{i+1}}\cdots x_{i_n}^{m_n}\ptl_{x_i}\;\mbox{is a root vector
corresponding to the
root}\;\sum_{r=i+1}^nm_r\ves_r-\ves_i.\eqno(3.10)$$ Expressions
(3.5) and (3.10) imply:\psp

{\bf Lemma 3.3}. {\it The subspace ${\cal L}_0({\cal T}^d)$ forms
an ideal of ${\cal L}_1({\cal T}^d)$ and is a direct sum of root
subspaces of dimension 1 with respect to $H$. Any abelian ideal of
${\cal L}_1({\cal T}^d)$ is a direct sum of some root subspaces of
${\cal L}_0({\cal T}^d)$}.\psp

By (3.5), the set of roots of ${\cal L}_0({\cal T}^d)$: $$R({\cal
L}_0({\cal T}^d))=\{-\ves_\es,\sum_{\iota_s\in{\cal D}_i}
j_s\ves_s-\ves_i\mid\iota_\es\in\Psi,\;
j_s\in\mbb{N},\;\sum_{\iota_s\in{\cal D}_i}j_s\kappa_{i,s}\leq
\kappa_i\}\eqno(3.11)$$

For convenience, we denote
$$v_{-\ves_\es}=\ptl_{x_\es},\;\;v_{_{\sum_{\iota_s\in{\cal D}_i}
j_s\ves_s-\ves_i}}=(\prod_{\iota_s\in{\cal D}_i} x^{j_s}_s
)\ptl_{x_i}.\eqno(3.12)$$ Suppose that ${\cal I}$ is a nonzero
abelian ideal of ${\cal L}_1({\cal T}^d)$. Define the set of roots
of ${\cal I}$ by:
$$R({\cal I})=\{\al\in R({\cal L}_0({\cal
T}^d))\mid v_\al\in{\cal I}\}.\eqno(3.13)$$ Then
$${\cal I}=\mbox{Span}\:\{v_\al\mid \al\in R({\cal
I})\}.\eqno(3.14)$$ Determining ${\cal I}$ is equivalent to
determining all possible $R({\cal I})$.

Set
$$\Phi=\{\iota_i\in{\cal N}\mid d({\cal
E}_{1,i})=\{1\}\}\bigcup\{\iota_1\}.\eqno(3.15)$$ Recall that a
subset $S$ of nodes is called  independent if each element of $S$
is not a descendent of any other element in $S$. For
$\iota_i\in{\cal N}$, we define
$$\Re_i=\{\{j_s\mid \iota_s\in{\cal D}_i\}\subset\mbb{N}\mid\sum_{\iota_s\in{\cal
D}_i}j_s\kappa_{i,s}\leq \kappa_i\}.\eqno(3.16)$$ For convenience,
we denote
$$\prod_{(\iota_r,\iota_s)\in{\cal
E}_{i,i}}d[(\iota_r,\iota_s)]=1.\eqno(3.17)$$ Moreover, we define
a partial ordering ``$\prec$'' on $\Re_i$ by
$$\vec j=\{j_s\mid \iota_s\in{\cal D}_i\}\prec \{l_s\mid \iota_s\in{\cal
D}_i\}\eqno(3.18)$$ if
$$\sum_{\iota_i\neq \iota_r\in{\cal C}_{i,p(m)}}j_r\prod_{(\iota_{\es},\iota_s)\in{\cal
E}_{r,p(m)}}d[(\iota_\es,\iota_s)]\leq \sum_{\iota_i\neq
\iota_r\in{\cal C}_{i,p(m)}}l_r\prod_{(\iota_\es,\iota_s)\in{\cal
E}_{r,p(m)}}d[(\iota_\es,\iota_s)]\eqno(3.19)$$ for $
\iota_m\in{\cal D}_i.$ Define
$$x^{\vec j}=\prod_{\iota_s\in{\cal D}_i} x^{j_s}_s\qquad\for\;\;\vec j=\{j_s\mid \iota_s\in{\cal
D}_i\}\in\Re_i.\eqno(3.20)$$ For convenience, we treat
$$\Re_\es=\{0\},\qquad x^0=1\qquad\for\;\;\iota_\es\in \Psi.\eqno(3.21)$$

 For any
independent subset $S$ of $\Phi$ (cf. (3.15)), we define
$${\cal I}(S)=\mbox{Span}\:\{x^{\vec
j}\ptl_{x_s}\mid\iota_i\in S,\;\vec j\in \Re_i,\;\iota_s\in{\cal
C}_i\}.\eqno(3.22)$$ It can be verified that ${\cal I}(S)$ is an
abelian ideal of ${\cal L}_1({\cal T}^d)$.
Set$$\Omega=\{\iota_i\in{\cal N}\mid (\iota_1,\iota_i)\in{\cal
E}\}.\eqno(3.23)$$ By similar arguments as those in the proofs of
Lemma 2.5 and Theorem 2.6, we have:\psp

{\bf Theorem 3.5}. {\it Any maximal abelian ideal of ${\cal
L}_1({\cal T}^d)$ is of the form ${\cal I}(S)$ for some maximal
independent subset $S$ of $\Phi$. Moreover, ${\cal L}_1({\cal
T}^d)$ has the unique maximal abelian ideal ${\cal I}(\Omega)$ if
$\Omega\bigcap\Phi=\emptyset$.}\psp

{\bf Example 3.1}. Suppose that ${\cal T}={\cal
T}_{E^{n_0}_{n_1,n_2}}$ (diagram before (1.4)). If $d({\cal
E})=\{1\}$, then $\Phi={\cal N}$ and the number of maximal abelian
ideals of ${\cal L}_1({\cal T}^d_{E^{n_0}_{n_1,n_2}})$ is
$n_0+n_1n_2$. Moreover,
$$\dim {\cal L}_0({\cal
T}^d_{E^{n_0}_{n_1,n_2}})=n_0(n_1+n_2)+\frac{n_0^2+n_1^2+n_2^2+n_0+n_1+n_2}{2}.\eqno(3.24)$$

When $n_0>1$ and $d[(\iota_1,\iota_2)]=m>1$,  ${\cal L}_1({\cal
T}^d_{E^{n_0}_{n_1,n_2}})$ has a unique maximal abelian ideal. If
we in addition  assume $d({\cal
E}\setminus\{(\iota_1,\iota_2)\})=\{1\}$, then
\begin{eqnarray*}\hspace{1.5cm}\dim {\cal L}_0({\cal
T}^d_{E^{n_0}_{n_1,n_2}})&=&{n_0+n_1+n_2+m-1\choose m}+
(n_0-1)(n_1+n_2)\\ & &
+\frac{n_0^2+n_1^2+n_2^2-n_0+n_1+n_2}{2}.\hspace{4.1cm}(3.25)\end{eqnarray*}
\vspace{0.1cm}

 Given $\iota_i\in{\cal N}$ and $\vec j\in \Re_i$,
we define
$${\cal I}[i,\vec j]=\mbox{Span}\:\{x^{\vec
l}\ptl_{x_r}\mid \iota_r\in{\cal C}_i,\;\vec j\succ\vec l\in
\Re_i\}.\eqno(3.26)$$ If $\iota_i\in\Phi$, ${\cal I}[i,\vec j]$ is
an abelian ideal of ${\cal L}_1({\cal T}^d)$ generated by $x^{\vec
j}\ptl_{x_i}$. A subset ${\cal K}$ of $\Re_i$ is called {\it
independent} if ${\cal K}$ itself is the set of maximal elements
in $\cal K$. Pick a non-empty independent set $S\subset \Phi$.
Define
$$W=\bigcup_{\iota_i\in S}{\cal C}_i.\eqno(3.27)$$
For each $\iota_i\in S$ and $\iota_\es\in{\cal C}_i$, we choose
independent subsets ${\cal K}_\es\subset \Re_i$ such that ${\cal
K}_i\neq\emptyset$ and
$$\vec j\not\prec\vec l\;\;\for\;\;\vec j\in{\cal K}_r,\;\vec
l\in{\cal K}_s\;\;\mbox{if}\;\;\iota_s\in{\cal
C}_i\;\mbox{and}\;\iota_s\neq\iota_r\in{\cal C}_s.\eqno(3.28)$$
Such an pair $(S,,\{{\cal K}_r\mid \iota_r\in W\})$ is called an
{\it admissible pair}.
 Moreover, we define
$${\cal I}[S,\{{\cal K}_r\mid \iota_r\in W\}]=\sum_{\iota_s\in
W}\sum_{\vec j\in {\cal K}_s}{\cal I}[s,\vec j].\eqno(3.29)$$ Then
${\cal I}[S,\{{\cal K}_r\mid \iota_r\in W\}]$ is a nonzero abelian
ideal of ${\cal L}_1({\cal T}^d)$ . By similar arguments as those
between (2.44)-(2.51), we obtain: \psp

{\bf Theorem 3.6}. {\it Any nonzero abelian ideal of ${\cal
L}_1({\cal T}^d)$ must be of the form ${\cal I}[S,\{{\cal K}_r\mid
\iota_r\in W\}]$ for some admissible pair $(S,\{{\cal K}_r\mid
\iota_r\in W\})$ constructed in the above paragraph. So there is a
one-to-one correspondence between the set of nonzero abelian
ideals of ${\cal L}_1({\cal T}^d)$ and the set of admissible
pairs.}\psp

\section{First-Order Differential Equations}

In this section, we will solve the equation (1.11) subject to the
condition (1.12).

First we want to solve the following evolution partial
differential equation:
$$u_t=(\ptl_{x_1}+x_1^{m_1}\ptl_{x_2}+x_2^{m_2}\ptl_{x_3}+\cdots+
x_{n-1}^{m_{n-1}}\ptl_{x_n})(u)\eqno(4.1)$$ subject to the
condition:
$$u(0,x_1,...,x_n)=f(x_1,...,x_n),\eqno(4.2)$$
where $f(x_1,x_2,...,x_n)$ is a smooth function.
 The equation (4.1) is related to
the Lie algebra $L_0({\cal T}_{A_n}^d)$ with
$$d[(\iota_i,\iota_{i+1})]=m_i\qquad\for\;\;i\in\ol{1,n-1}.\eqno(4.3)$$

For convenience, we denote $m_0=1$ and $x_0=1$. Set
$$D_i=t\sum_{r=0}^{i-1}x_r^{m_r}\ptl_{x_{r+1}}\qquad\for\;\;i\in\ol{1,n}.\eqno(4.4)$$
In order to solve the above initial value problem, we need to
factorize
$$e^{D_n}=e^{\eta_n(t)\ptl_{x_n}}e^{\eta_{n-1}(t)\ptl_{x_{n-1}}}\cdots
e^{\eta_1(t)\ptl_{x_1}},\eqno(4.5)$$ where
$\eta_1(t),...,\eta_n(t)$ are suitable functions to be
constructed. This operator applied to the initial condition (4.2)
will yield the solution of the initial value problem.

Denote
$$ A=D_n,\;\;B=-tx_{n-1}^{m_{n-1}}\ptl_{x_n}.\eqno(4.6)$$
Thus
$$D_{n-1}=D_n+B=A+B.\eqno(4.7)$$
 The Lie ideal of $L_0({\cal T}_{A_n}^d)$ generated by $B$
 (with $0\neq t\in\mbb{R}$)
is the abelian ideal: \begin{eqnarray*}\hspace{1cm}{\cal
I}(\iota_n)&=&\mbox{Span}\:\{x_1^{j_1}x_2^{j_2}\cdots
x_{n-1}^{j_{n_1}}\ptl_{x_n}\mid
(j_1,...,j_{n-1})\in\mbb{N}^{\:n-1},\\ &
&\sum_{i=1}^{n-1}j_im_0m_1\cdots m_{i-1}\leq m_0m_1\cdots
m_{n-1}\}.\hspace{5.2cm}(4.8)\end{eqnarray*} According to the
Campbell-Hausdorff formula (e.g., cf. Section 5 of Chapter 5 in
[J2]), (4.8) being an abelian ideal implies
$$\ln e^Ae^B=A+B+\sum_{r=1}^{\infty}a_r(\ad A)^r(B),\qquad
a_r\in\mbb{R}.\eqno(4.9)$$ On the other hand,
$$\ln
e^Ae^B=\sum_{m=1}^{\infty}\sum_{p_i,q_i\in\mbb{N},\;p_i+q_i>0,\;i\in\ol{1,m}}\frac{(-1)^{m-1}}{m}
\frac{A^{p_1}}{p_1!}\frac{B^{q_1}}{q_1!}\cdots\frac{A^{p_m}}{p_m!}\frac{B^{q_m}}{q_m!}.
\eqno(4.10)$$  Note that
$$(\ad
A)^r(B)=\sum_{i=0}^r(-1)^{r-i}{r\choose{i}}A^iBA^{r-i}.\eqno(4.11)$$

For any positive integer $k$, the coefficients of $A^kB$ in (4.9)
and (4.10) imply
$$a_k=\sum_{m=1}^{\infty}\;\sum_{p_1,...,p_m\in\mbb{N},p_1+\cdots+p_m=k+1-m}\frac{(-1)^{m-1}}{m}
\frac{1}{(p_1+1)!\cdots (p_{m-1}+1)!p_m!}.\eqno(4.12)$$ Set
$$b_m=\sum_{p_1,...,p_m\in\mbb{N},\;p_1+\cdots+p_m=k+1-m}
\frac{k!}{(p_1+1)!\cdots (p_{m-1}+1)!p_m!}.\eqno(4.13)$$ Then
$$b_m=m^k-\sum_{i=1}^{m-1}{m-1\choose{i}}b_{m-i}.\eqno(4.14)$$
$$b_1=1=(y\ptl_y)^k(y)|_{y=1},\;\;b_2=2^k-1=(y\ptl_y)^k(y^2-y)|_{y=1}=(y\ptl_y)^k[y(y-1)]|_{y=1},
\eqno(4.15)$$
$$b_3=3^k-2b_2-1=3^k-2(2^k-1)-1=3^k-2\cdot
2^k+1=(y\ptl_y)^k[y(y-1)^2]|_{y=1}.\eqno(4.16)$$ Suppose
$$b_i=(y\ptl_y)^k[y(y-1)^{i-1}]|_{y=1}\qquad\for\;\;i<r.\eqno(4.17)$$

Now
\begin{eqnarray*}\hspace{1cm}b_r&=&r^k-\sum_{i=1}^{r-1}{r-1\choose{i}}b_{r-i}=
(t\ptl_t)^k[y^r-\sum_{i=1}^{r-1}{r-1\choose i}y(y-1)^{r-i-1}]|_{y=1}\\
&=&(y\ptl_y)^k[y^r-y(\sum_{i=0}^{r-1}{r-1\choose
i}(y-1)^{r-i-1}-(y-1)^{r-1})]|_{y=1}
\\ &=&(y\ptl_y)^k[y^r-y((y-1+1)^{r-1}-(y-1)^{r-1})]|_{y=1}\\
&=&(y\ptl_y)^k[y^r-y(y^{r-1}-(y-1)^{r-1})]|_{y=1}\\ &
=&(y\ptl_y)^k[y(y-1)^{r-1}]|_{y=1}. \hspace{8.4cm}(4.18)
\end{eqnarray*}
Thus
$$b_m=(t\ptl_t)^k[y(y-1)^{m-1}]|_{y=1}\qquad\for\;\;m\in
\mbb{N}+1\eqno(4.19)$$ by induction.

Observe \begin{eqnarray*}a_k&=&
\frac{1}{k!}\sum_{m=1}^{\infty}\frac{(-1)^{m-1}}{m}b_m=
\frac{1}{k!}\sum_{m=1}^{\infty}\frac{(-1)^{m-1}}{m}(y\ptl_y)^k[y(y-1)^{m-1}]|_{y=1}\\
&=&\left[\frac{(y\ptl_y)^k}{k!}\left(\frac{y}{y-1}\sum_{m=1}^{\infty}\frac{(-1)^{m-1}}{m}
(y-1)^m\right)\right]|_{y=1}
\\ &=&\left[\frac{(y\ptl_y)^k}{k!}\left(\frac{y\ln(1+y-1)}{y-1}\right)\right]|_{y=1}
= \left[\frac{(y\ptl_y)^k}{k!}\left(\frac{y\ln y
}{y-1}\right)\right]|_{y=1}.\hspace{2.9cm}(4.20)\end{eqnarray*}
Denote $a_0=1$. Hence
\begin{eqnarray*}\hspace{1cm}\sum_{k=0}^{\infty}a_kx^k&=&\sum_{k=0}^{\infty}
\left[\frac{(xy\ptl_y)^k}{k!}\left(\frac{y\ln
y}{y-1}\right)\right]|_{y=1}= e^{xy\ptl_y}\left(\frac{y\ln
y}{y-1}\right)|_{y=1}\\ &=&\left(\frac{e^xy\ln
e^xy}{e^xy-1}\right)|_{y=1}=\frac{e^x\ln
e^x}{e^x-1}=\frac{x}{1-e^{-x}}.\hspace{4.3cm}(4.21)\end{eqnarray*}

Since we have
$$e^Ae^B=e^{A+\sum_{i=0}^{\infty}a_i(\ad A)^i(B)},\eqno(4.22)$$
we have $$e^{D_{n-1}-B}=e^A=e^{A+\sum_{i=0}^{\infty}a_i(\ad
A)^i(B)}e^{-B}=e^{D_{n-1}+\sum_{i=1}^{\infty}a_i(\ad
A)^i(B)}e^{-B}.\eqno(4.23)$$ Note
$$\ad A|_{{\cal I}(\iota_n)}=\ad D_{n-1}|_{{\cal I}(\iota_n)}=\ad(D_{n-1}+
\sum_{i=1}^{\infty}a_i(\ad A)^i(B))|_{{\cal
I}(\iota_n)}.\eqno(4.24)$$ Viewing $\sum_{i=1}^{\infty}a_i(\ad
A)^i(B)$ as $-B$, we have:
$$e^{D_{n-1}+\sum_{i=1}^{\infty}a_i(\ad
A)^i(B)}=e^{D_{n-1}-\sum_{i,j=1}^{\infty}a_ia_j(\ad
A)^{i+j}(B)}e^{\sum_{i=1}^{\infty}a_i(\ad A)^i(B)}\eqno(4.25)$$ by
(4.21). Thus
$$e^A=e^{D_{n-1}-\sum_{i,j=1}^{\infty}a_ia_j(\ad
A)^{i+j}(B)}e^{-B+\sum_{i=1}^{\infty}a_i(\ad
A)^i(B)}.\eqno(4.26)$$ By induction, we obtain:
\begin{eqnarray*}\hspace{1cm}e^A&=&e^{D_{n-1}-(-1)^m\sum_{i_1,i_2,...,i_m=1}^{\infty}a_{i_1}
\cdots a_{i_m}(\ad A)^{i_1+\cdots+i_m}(B)}\\ & &\times
e^{-B+\sum_{r=0}^{m-1}(-1)^{r-1}\sum_{i_1,..,i_r=1}^{\infty}a_{i_1}\cdots
a_{i_r} (\ad
A)^{i_1+\cdots+i_r}(B)}.\hspace{4cm}(4.27)\end{eqnarray*} Since
$L_0({\cal T}_{A_n}^d)$ is a nilpotent Lie algebra, we get
$$e^{D_n}=e^A=e^{D_{n-1}}e^{-B+\sum_{r=1}^{\infty}(-1)^r\sum_{i_1,..,i_r=1}^{\infty}a_{i_1}\cdots
a_{i_r} (\ad A)^{i_1+\cdots+i_r}(-B)}.\eqno(4.28)$$

Note
\begin{eqnarray*} \hspace{3cm}& &1+\sum_{r=1}^{\infty}(-1)^r\sum_{i_1,..,i_r=1}^{\infty}a_{i_1}\cdots
a_{i_r}x^{i_1+\cdots+i_r}\\
&=&1+\sum_{r=1}^{\infty}(-1)^r(\sum_{i=1}^{\infty}a_ix^i)^r=\frac{1}{1+\sum_{i=1}^{\infty}a_ix^i}
\\
&=&\frac{1}{x/(1-e^{-x})}=\frac{1-e^{-x}}{x}\hspace{6.8cm}(4.29)\end{eqnarray*}
by (4.19). Denote
$$\vt(x)=\frac{1-e^{-x}}{x}=\int_{-1}^0e^{yx}dy=\sum_{i=1}^{\infty}
\frac{(-1)^{i-1}}{i!}x^{i-1}.\eqno(4.30)$$ According to (4.24),
$$(\ad A)^i(-B)=(\ad
D_{n-1})^i(x_{n-1}^{m_{n-1}}\ptl_{x_n})=D_{n-1}^i(x_{n-1}^{m_{n-1}})\ptl_{x_n}.\eqno(4.31)$$
Now (4.28)-(4.31) imply:\psp

{\bf Lemma 4.1}. {\it We have the following operator formula:}
$$e^{D_n}=e^{D_{n-1}}e^{t\vt(D_{n-1})(x_{n-1}^{m_{n-1}})\ptl_{x_n}}.\eqno(4.32)$$
\vspace{0.1cm}

We remark that the above formula holds for any $D_i$. Next we want
to find $\eta_1(t),...,\eta_n(t)$ in (4.5). We set
$$\xi_1(t)=t,\;\;\xi_i(t)=t\vt(D_{i-1})(x_{i-1}^{m_{i-1}})\qquad\for\;\;i\in\ol{2,n}.
\eqno(4.33)$$ By (4.4) and (4.31), we get
$$e^{D_i}=e^{\xi_1(t)\ptl_{x_1}}e^{\xi_2(t)\ptl_{x_2}}\cdots
e^{\xi_i(t)\ptl_{x_i}}\qquad\for\;\;i\in\ol{1,n}.\eqno(4.34)$$
Moreover, we define
$$\eta_1(t)=t,\;\;\eta_i(t)=e^{D_{i-1}}(\xi_i(t)).\eqno(4.35)$$

For a smooth function $g(x)$, Taylor's Theorem says
$$e^{z\ptl_x}(g(x))=g(x+z).\eqno(4.36)$$
Moreover,
$$e^ABe^{-A}=e^{\ad A}(B)\qquad\for\;\;A,B\in\mbb{A}.\eqno(4.37)$$
Since $D_1=t\ptl_{x_1}$, we have
$$\xi_2(t)=t\int_{-1}^0e^{ty_1\ptl_{x_1}}(x_1^{m_1})dy_1=t\int_{-1}^0(x_1+ty_1)^{m_1}dy_1=
\int_{-t}^0(x_1+y_1)^{m_1}dy_1.\eqno(4.38)$$ Furthermore,
\begin{eqnarray*}\hspace{2cm}\eta_2&=&e^{D_1}(\xi_2)=e^{t\ptl_{x_1}}\left[\int_{-t}^0(x_1+y_1)^{m_1}dy_1\right]
\\
&=&\int_{-t}^0(x_1+t+y_1)^{m_1}dy_1=\int_0^t(x_1+y_1)^{m_1}dy_1\hspace{4cm}(4.39)
\end{eqnarray*}
by (4.36). Now,
\begin{eqnarray*}\xi_3(t)&=&t\int_{-1}^0e^{y_2D_2}(x_2^{m_2})dy_2
=t\int_{-1}^0e^{\xi_1(ty_2)\ptl_{x_1}}e^{\xi_2(ty_2)\ptl_{x_2}}(x_2^{m_2})dy_2\\
&=&\int_{-t}^0[e^{y_2\ptl_{x_1}}e^{\xi_2(y_2)\ptl_{x_2}}e^{-y_2\ptl_{x_1}}]
e^{y_2\ptl_{x_1}}(x_2^{m_2})dy_2=\int_{-t}^0e^{\ad
y_2\ptl_{x_1}}(e^{\xi_2(y_2)\ptl_{x_2}})(x_2^{m_2})dy_2\\
&=&\int_{-t}^0e^{e^{\ad
y_2\ptl_{x_1}}(\xi_2(y_2)\ptl_{x_2})}(x_2^{m_2})dy_2=\int_{-t}^0e^{e^{
y_2\ptl_{x_1}}(\xi_2(y_2))\ptl_{x_2}}(x_2^{m_2})dy_2\\
&=&\int_{-t}^0e^{\eta_2(y_2)\ptl_{x_2}}(x_2^{m_2})dy_2=\int_{-t}^0(x_2+\eta_2(y_2))^{m_2}dy_2\\
&=&\int_{-t}^0(x_2+\int_0^{y_2}(x_1+y_1)^{m_1}dy_1)^{m_2}dy_2\hspace{7cm}(4.40)\end{eqnarray*}
by (4.35) and (4.37). We calculate
\begin{eqnarray*}\eta_3(t)&=&e^{D_2}(\xi_3(t))=e^{t\ptl_{x_1}}e^{\xi_2(t)\ptl_{x_2}}
[\int_{-t}^0(x_2+\int_0^{y_2}(x_1+y_1)^{m_1}dy_1)^{m_2}dy_2]\\
&=&[e^{t\ptl_{x_1}}e^{\xi_2(t)\ptl_{x_2}}e^{-t\ptl_{x_1}}]e^{t\ptl_{x_1}}
[\int_{-t}^0(x_2+\int_0^{y_2}(x_1+y_1)^{m_1}dy_1)^{m_2}dy_2]\\
&=&e^{\eta_2(t)\ptl_{x_2}}[\int_{-t}^0(x_2+\int_0^{y_2}(x_1+t+y_1)^{m_1})^{m_2}dy_2]\\
&=&\int_{-t}^0(x_2+\eta_2(t)+\int_t^{y_2+t}(x_1+y_1)^{m_1}dy_1)^{m_2}dy_2\\
&=&\int_{-t}^0(x_2+\int_0^t(x_1+y_1)^{m_1}dy_1+\int_t^{y_2+t}(x_1+y_1)^{m_1})^{m_2}dy_2\\
&=& \int_{-t}^0(x_2+\int_0^{y_2+t}(x_1+y_1)^{m_1}dy_1)^{m_2}dy_2\\
&
=&\int_0^t(x_2+\int_0^{y_2}(x_1+y_1)^{m_1}dy_1)^{m_2}dy_2.\hspace{6.8cm}(4.41)\end{eqnarray*}

Suppose
$$\xi_i(t)=\int_{-t}^0(x_{i-1}+\int^{y_{i-1}}_0(x_{i-2}+...+\int_0^{y_2}(x_1+y_1)^{m_1}dy_1...)^{m_{i-2}}dy_{i-2})^{m_{i-1}
}dy_{i-1}\eqno(4.42)$$ and
$$\eta_i(t)=\int_0^t(x_{i-1}+\int^{y_{i-1}}_0(x_{i-2}+...+\int_0^{y_2}(x_1+y_1)^{m_1}dy_1...)^{m_{i-2}}dy_{i-2})^{m_{i-1}
}dy_{i-1}\eqno(4.43)$$ for $2\leq i\leq r<n$. Note
$$e^{D_r}=e^{D_{r-1}}e^{\xi_r(t)\ptl_{x_r}}=[e^{D_{r-1}}e^{\xi_r(t)\ptl_{x_r}}e^{-D_{r-1}}]
e^{D_{r-1}}=e^{\eta_r(t)\ptl_{x_r}}e^{D_{r-1}}\eqno(4.44)$$ by
(4.30) and (4.35). By induction, we have:
$$e^{D_r}=e^{\eta_r(t)\ptl_{x_r}}e^{\eta_{r-1}(t)\ptl_{x_{r-1}}}\cdots
e^{\eta_1(t)\ptl_{x_1}}.\eqno(4.45)$$ Hence
\begin{eqnarray*}
\xi_{r+1}(t)&=&t\int_{-1}^0e^{y_rD_r}(x_r^{m_r})dy_r
=\int_{-1}^0e^{\eta_r(ty_r)\ptl_{x_r}}e^{\eta_{r-1}(ty_r)\ptl_{x_{r-1}}}\cdots
e^{\eta_1(ty_r)\ptl_{x_1}}(x_r^{m_r})tdy_r\\
&=&\int_{-t}^0e^{\eta_r(y_r)\ptl_{x_r}}(x_r^{m_r})dy_r
=\int_{-t}^0(x_r+\eta_r(y_r))^{m_r}dy_r\\
&=&\int_{-t}^0(x_r+\int_0^{y_r}(x_{r-1}+\int^{y_{r-1}}_0
(x_{r-2}+...+\int_0^{y_2}(x_1+y_1)^{m_1}\\ & &
dy_1...)^{m_{r-2}}dy_{r-2})^{m_{r-1}
}dy_{r-1})^{m_r}dy_r.\hspace{6.7cm}(4.46)\end{eqnarray*} By
induction, (4.42) holds for any $2\leq i<n$. Furthermore, using
(4.44) and applying (4.36) repeatedly, we obtain
\begin{eqnarray*}& &
\eta_{r+1}(t)\\
&=&e^{D_r}(\xi_{r+1}(t))=e^{\eta_r(t)\ptl_{x_r}}e^{\eta_{r-1}(t)\ptl_{x_{r-1}}}\cdots
e^{\eta_1(t)\ptl_{x_1}}(\xi_{r+1}(t))\\
&=&e^{\eta_r(t)\ptl_{x_r}}e^{\eta_{r-1}(t)\ptl_{x_{r-1}}}\cdots
e^{t\ptl_{x_1}}(\int_{-t}^0(x_r
+\int_0^{y_r}(x_{r-1}+\int^{y_{r-1}}_0 (x_{r-2}+...\\ &
&+\int_0^{y_2}(x_1+y_1)^{m_1}dy_1...)^{m_{r-2}}dy_{r-2})^{m_{r-1}
}dy_{r-1})^{m_r}dy_r)\\ &=&
e^{\eta_r(t)\ptl_{x_r}}e^{\eta_{r-1}(t)\ptl_{x_{r-1}}}\cdots
e^{\eta_2(t)\ptl_{x_2}} (\int_{-t}^0(x_r
+\int_0^{y_r}(x_{r-1}+\int^{y_{r-1}}_0
(x_{r-2}+...\\
&
&+\int_0^{y_2}(x_1+t+y_1)^{m_1}dy_1...)^{m_{r-2}}dy_{r-2})^{m_{r-1}
}dy_{r-1})^{m_r}dy_r)\\
&=&e^{\eta_r(t)\ptl_{x_r}}e^{\eta_{r-1}(t)\ptl_{x_{r-1}}}\cdots
e^{\eta_3(t)\ptl_{x_3}}
(\int_{-t}^0(x_r+\int_0^{y_r}(x_{r-1}+\int^{y_{r-1}}_0
(x_{r-2}+...\\
& &
+\int_0^{y_3}(x_2+\eta_2(t)+\int_t^{y_2+t}(x_1+y_1)^{m_1}dy_1)^{m_2}dy_2...)^{m_{r-2}}dy_{r-2})^{m_{r-1}
}dy_{r-1})^{m_r}dy_r)\\
&=&e^{\eta_r(t)\ptl_{x_r}}e^{\eta_{r-1}(t)\ptl_{x_{r-1}}}\cdots
e^{\eta_3(t)\ptl_{x_3}}(\int_{-t}^0(x_r
+\int_0^{y_r}(x_{r-1}+\int^{y_{r-1}}_0
(x_{r-2}+...\\
&
&+\int_0^{y_3}(x_2+\int_0^{y_2+t}(x_1+y_1)^{m_1}dy_1)^{m_2}dy_2...)^{m_{r-2}}dy_{r-2})^{m_{r-1}
}dy_{r-1})^{m_r}dy_r)\\
&=&e^{\eta_r(t)\ptl_{x_r}}e^{\eta_{r-1}(t)\ptl_{x_{r-1}}}\cdots
e^{\eta_3(t)\ptl_{x_3}}
(\int_{-t}^0(x_r+\int_0^{y_r}(x_{r-1}+\int^{y_{r-1}}_0
(x_{r-2}+...\\
& &
+\int_t^{y_3+t}(x_2+\int_0^{y_2}(x_1+y_1)^{m_1}dy_1)^{m_2}dy_2...)^{m_{r-2}}dy_{r-2})^{m_{r-1}
}dy_{r-1})^{m_r}dy_r)=\cdots\\ &=&e^{\eta_r(t)\ptl_{x_r}}
(\int_{-t}^0(x_r+\int_t^{y_r+t}(x_{r-1}+\int^{y_{r-1}}_0
(x_{r-2}+...\\ & &
+\int_t^{y_3}(x_2+\int_0^{y_2}(x_1+y_1)^{m_1}dy_1)^{m_2}dy_2...)^{m_{r-2}}dy_{r-2})^{m_{r-1}
}dy_{r-1})^{m_r}dy_r)\\
&=&\int_{-t}^0(x_r+\eta_r(t)+\int_t^{y_r+t}(x_{r-1}+\int^{y_{r-1}}_0
(x_{r-2}+...+\int_0^{y_2}(x_1+y_1)^{m_1}\\ & &
dy_1)^{m_2}dy_2...)^{m_{r-2}}dy_{r-2})^{m_{r-1}
}dy_{r-1})^{m_r}dy_r\\ &=&
\int_{-t}^0(x_r+\int_0^{y_r+t}(x_{r-1}+\int^{y_{r-1}}_0
(x_{r-2}+...+\int_0^{y_2}(x_1+y_1)^{m_1}\\ & &
dy_1)^{m_2}dy_2...)^{m_{r-2}}dy_{r-2})^{m_{r-1}
}dy_{r-1})^{m_r}dy_r\hspace{9cm}\end{eqnarray*}
\begin{eqnarray*}
&=&\int_0^t(x_r+\int_0^{y_r}(x_{r-1}+\int^{y_{r-1}}_0
(x_{r-2}+...+\int_0^{y_2}(x_1+y_1)^{m_1}\\ & &
dy_1)^{m_2}dy_2...)^{m_{r-2}}dy_{r-2})^{m_{r-1}
}dy_{r-1})^{m_r}dy_r.\hspace{6.7cm}(4.47)\end{eqnarray*} This
shows that (4.43) holds for any $2\leq i\leq n$. Thus (4.5) is
determined.

Applying to the initial condition (4.2), we get the solution
\begin{eqnarray*}\hspace{1cm}u&=&e^{D_n}(f(x_1,x_2,...,x_n)=e^{\eta_n(t)\ptl_{x_n}}\cdots
e^{\eta_1(t)\ptl_{x_1}}(f(x_1,x_2,...,x_n))\\ &=&
f(x_1+\eta_1(t),x_2+\eta_2(t),...,x_n+\eta_n(t)).,\hspace{5.6cm}(4.48)\end{eqnarray*}
because $u_t=\ptl_t
e^{D_n}(f)=(\sum_{i=0}^{n-1}x_i^{m_i}\ptl_{x_{i+1}})e^{D_n}(f)=
(\sum_{i=0}^{n-1}x_i^{m_i}\ptl_{x_{i+1}})(u)$. Since (4.1) is a
first-order partial differential equation, we have:\psp

{\bf Proposition 4.2}. {\it Let $f(x_1,x_2,...,x_n)$ be a first
order differentiable function. Then
$$u=f(x_1+\eta_1(t),x_2+\eta_2(t),...,x_n+\eta_n(t))\eqno(4.49)$$
is the solution of the equation (4.1) subject to (4.48), where
$\eta_1=t$ and $\eta_i(t)$ are given in (4.2) for $i\in\ol{2,n}$.}
\psp

The above process of solving the initial value problem (4.1) and
(4.2) is necessary in order to solve the following more general
problem and the initial value problem (1.9) and (1.10) in next
section. Let ${\cal T}^d$ be any tree diagram with $n$ nodes.
Consider the equation:
$$u_t=(\ptl_{x_1}+\sum_{(\iota_i,\iota_j)\in{\cal
E}}x_i^{d[(\iota_i,\iota_j)]}\ptl_{x_j})(u)\eqno(4.50)$$ subject
to the condition (4.2). Take $\eta_1(t)=t$ and redefine
$\{\eta_i(t)\mid i\in\ol{2,n}\}$ as follows. Suppose ${\cal
C}_i=\{\iota_{i_1},\iota_{i_2},...,\iota_r\}$. We define
\begin{eqnarray*}\hspace{2cm}\xi_i(t)&=&\int_{-t}^0(x_{i_{r-1}}+\int^{y_{i_{r-1}}}_0(x_{i_{r-2}}+...+\int_0^{y_{i_2}}
(x_1+y_1)^{d[(\iota_1,\iota_2)]}\\ &&dy_1...
)^{d[(\iota_{i_{r-2}},\iota_{i_{r-1}})]}dy_{i_{r-2}})^{d[(\iota_{i_{r-1}},\iota_{i_r})]}
dy_{i_{r-1}}\hspace{4.2cm}(4.51)\end{eqnarray*} and
\begin{eqnarray*}\hspace{2cm}\eta_i(t)&=&\int_0^t(x_{i_{r-1}}+\int^{y_{i_{r-1}}}_0(x_{i_{r-2}}+...+\int_0^{y_{i_2}}
(x_1+y_1)^{d[(\iota_1,\iota_2)]}\\ &&dy_1...
)^{d[(\iota_{i_{r-2}},\iota_{i_{r-1}})]}dy_{i_{r-2}})^{d[(\iota_{i_{r-1}},\iota_{i_r})]}
dy_{i_{r-1}}.\hspace{4cm}(4.52)\end{eqnarray*}

Note
$$\ptl_{x_r}(\xi_i(t))\neq 0\llra \iota_i\neq\iota_r\in{\cal
C}_i.\eqno(4.53)$$ Moreover,
$${\cal C}_i\subset\{\iota_1,\iota_2,...,\iota_i\}.\eqno(4.54)$$
 Set
$$D(t)=t\sum_{(\iota_i,\iota_j)\in{\cal
E}}x_i^{d[(\iota_i,\iota_j)]}\ptl_{x_j}.\eqno(4.55)$$ By
(4.32)-(4.34), (4.42) and (4.51)-(4.55), we obtain
$$ e^{D(t)}=e^{\xi_1(t)\ptl_{x_1}}e^{\xi_2(t)\ptl_{x_2}}\cdots
e^{\xi_n(t)\ptl_{x_n}}.\eqno(4.56)$$ Furthermore, (4.44), (4.53)
and (4.54) imply
$$e^{D(t)}=e^{\eta_n(t)\ptl_{x_n}}e^{\eta_{n-1}(t)\ptl_{x_{n-1}}}\cdots
e^{\xi_1(t)\ptl_{x_1}}.\eqno(4.57)$$ Therefore, we get our main
theorem in this section.\psp

{\bf Theorem 4.3}. {\it Let ${\cal T}^d$ be any tree diagram with
$n$ nodes and let $f(x_1,x_2,...,x_n)$ be a first-order
differentiable function. Then
$$u=f(x_1+\eta_1(t),x_2+\eta_2(t),...,x_n+\eta_n(t))\eqno(4.58)$$
is the solution of the equation (4.50) subject to (4.2), where
$\eta_1=t$ and $\eta_i(t)$ are given in (4.52) for
$i\in\ol{2,n}$.}\psp

{\bf Remark 4.4}. Since $d:{\cal E}\rta \mbb{N}+1$ is an arbitrary
map, we can solve more general problem of replacing exponential
functions by any first-order differentiable functions. Let
$\{g_i(x)\mid \iota_i\in{\cal N}\setminus \Psi\}$ be a set of
first-order differentiable functions. Consider the equation:
$$u_t=(\ptl_{x_1}+\sum_{(\iota_i,\iota_j)\in{\cal
E}}g_i(x_i)\ptl_{x_j})(u)\eqno(4.59)$$ subject to the condition
(4.2).

Suppose ${\cal C}_i=\{\iota_{i_1},\iota_{i_2},...,\iota_r\}$. We
define $$\eta_i(t)=\int_0^tg_{i_{r-1}}(x_{i_{r-1}}+
\int^{y_{i_{r-1}}}_0g_{i_{r-2}}(x_{i_{r-2}}+...+\int_0^{y_{i_2}}g_{i_1}
(x_1+y_1)dy_1...)dy_{i_{r-2}}) dy_{i_{r-1}}\eqno(4.60)$$ where $
i_1=1$ and $i_r=i$. Then (4.58) is the solution of the equation
(4.59) subject to (4.2), where $\eta_1=t$ and $\eta_i(t)$ are
given in (4.60) for $i\in\ol{2,n}$.

\section{High-Order Differential Equations}

In this section, we will solve the equation (1.9) subject to the
condition (1.10).

 Given a continuous function $f(x_1,x_2,...,x_n)$
on the region:
$$-a_i\leq x_i\leq a_i,\qquad
0<a_i\in\mbb{R},\qquad\for\;\;i\in\ol{1,n}.\eqno(5.1)$$ First we
want to solve the differential equation:
$$u_t=(\ptl_{x_1}^{m_1}+x_1\ptl_{x_2}^{m_2}+x_2\ptl_{x_3}^{m_3}
+\cdots+x_{n-1}\ptl_{x_n}^{m_n})(u)\eqno(5.2)$$ subject to the
initial condition:
$$u(0,x_1,...,x_n)=f(x_1,x_2,...,x_n)\qquad\for\;\;x_i\in[-a_i,a_i].\eqno(5.3)$$
Denote
$$D(t)=t(\ptl_{x_1}^{m_1}+x_1\ptl_{x_2}^{m_2}+x_2\ptl_{x_3}^{m_3}
+\cdots+x_{n-1}\ptl_{x_n}^{m_n}).\eqno(5.4)$$ Since
$$[\ptl_{x_i},x_j]=[-x_j,\ptl_{x_i}]=\dlt_{i,j},\eqno(5.5)$$
the linear map $\ast$ determined by
$$(x_i^l\ptl_{x_j}^k)^\ast=x_j^k(-\ptl_{x_i})^l\qquad\for\;\;i,j\in\ol{1,n},\;l,k\in\mbb{N}
\eqno(5.6)$$ is a Lie algebra automorphism of
$(\mbb{A},[\cdot,\cdot])$. We may write
$$z_i=\ptl_{x_i},\qquad
\ptl_{z_i}=-x_i\qquad\for\;\;i\in\ol{1,n}.\eqno(5.7)$$

Now
$$D(t)=t(z_1^{m-1}-z_2^{m_2}\ptl_{z_1}-z_3^{m_3}\ptl_{z_2}-\cdots
-z_{n-1}^{m_{n-1}}\ptl_{z_{n-2}}-z_n^{m_n}\ptl_{z_{n-1}}).\eqno(5.8)$$
Changing variables:
$$\zeta_i=\frac{z_i}{z_n^{m_{i+1}m_{i+2}\cdots
m_n}}\qquad\for\;\;i\in\ol{1,n-1},\eqno(5.9)$$ we have:
$$\ptl_{\zeta_i}=z_n^{m_{i+1}m_{i+2}\cdots
m_n}\ptl_{z_i}\qquad\for\;\;i\in\ol{1,n-1}\eqno(5.10)$$ and so
$$D(t)=-t(\zeta_1^{m_1}(-z_n^{m_1m_2\cdots
m_n})+\zeta_2^{m_2}\ptl_{\zeta_1}+\cdots+\zeta_{n-1}^{m_{n-1}}\ptl_{\zeta_{n-2}}
+\ptl_{\zeta_{n-1}}).\eqno(5.11)$$ Denote
$$D_i=-t(\zeta_i^{m_i}\ptl_{\zeta_{i-1}}+\zeta_{i+1}^{m_{i+1}}\ptl_{\zeta_i}+\cdots+\zeta_{n-1}^{m_{n-1}}\ptl_{\zeta_{n-2}}
+\ptl_{\zeta_{n-1}})\eqno(5.12)$$ for $i\in\ol{2,n-1}$. By Lemma
4.1,
$$e^{D(t)}=e^{D_2}e^{\vt(D_2)(\zeta_1^{m_1})tz_n^{m_1m_2\cdots
m_n}}.\eqno(5.13)$$

By (4.38)-(4.44), we have \begin{eqnarray*}
\xi_1&=&\vt(D_2)(\zeta_1^{m_1})tz_n^{m_1m_2\cdots
m_n}=-\int_t^0z_n^{m_1m_2\cdots
m_n}(\zeta_1+\int_0^{y_1}(\zeta_2+...\\ & &
+\int_0^{y_{n-2}}(\zeta_{n-1}+y_{n-1})^{m_{n-1}}dy_{n-1}...)
^{m_2}dy_2)^{m_1}dy_1\\
&=&\int_0^t(z_1+\int_0^{y_1}(z_2+...+\int_0^{y_{n-2}}(z_{n-1}+z^{m_n}y_{n-1})^{m_{n-1}}dy_{n-1}...)
^{m_2}dy_2)^{m_1}dy_1.\hspace{0.6cm}(5.14)\end{eqnarray*}
Similarly,
$$e^{D_i}=e^{D_{i+1}}e^{-t\vt(D_{i+1})(\zeta_i^{m_i})\ptl_{\zeta_{i-1}}}\qquad\for\;
\;i\in\ol{2,n-1}.\eqno(5.15)$$ We have
\begin{eqnarray*}\hspace{1cm}
\xi_i&=&-t\vt(D_{i+1})(\zeta_i^{m_i})\ptl_{\zeta_{i-1}}=
\int_t^0(\zeta_i+\int_0^{y_i}(\zeta_{i+1}+...\\ & &
+\int_0^{y_{n-2}}(\zeta_{n-1}+y_{n-1})^{m_{n-1}}dy_{n-1}...)
^{m_{i+1}}dy_{i+1})^{m_i}dy_i\ptl_{\zeta_{i-1}}\\
&=&\int_t^0(\zeta_i+\int_0^{y_i}(\zeta_{i+1}+...+\int_0^{y_{n-2}}(\zeta_{n-1}+y_{n-1})^{m_{n-1}}\\
&& dy_{n-1}...) ^{m_{i+1}}dy_{i+1})^{m_i}dy_iz_n^{m_im_{i+1}\cdots
m_n}
\ptl_{z_{i-1}}\\
&=&\int_t^0(z_i+\int_0^{y_i}(z_{i+1}+...+\int_0^{y_{n-2}}(z_{n-1}+z_n^{m_n}y_{n-1})^{m_{n-1}}\\
& & dy_{n-1}...)
^{m_{i+1}}dy_{i+1})^{m_i}dy_i\ptl_{z_{i-1}}.\hspace{7.4cm}(5.16)\end{eqnarray*}
We take
$$\xi_n=-t\ptl_{\zeta_n}=-tz_n^{m_n}\ptl_{z_{n-1}}.\eqno(5.17)$$

According to (5.13) and (5.15), we have
$$e^{D(t)}=e^{\xi_n}e^{\xi_{n-1}}\cdots
e^{\xi_2}e^{\xi_1}.\eqno(5.18)$$ Substituting (5.7) into
(5.15)-(5.17), we obtain:
\begin{eqnarray*} \hspace{1cm}& &\xi_1(t,\ptl_{x_1},...,\ptl_{x_n})=
\int_0^t(\ptl_{x_1}+\int_0^{y_1}(\ptl_{x_2}+...+\int_0^{y_{n-2}}(\ptl_{x_{n-1}}\\
&& +y_{n-1} \ptl_{x_n}^{m_n})^{m_{n-1}}dy_{n-1}...)
^{m_2}dy_2)^{m_1}dy_1,\hspace{6.6cm}(5.19)\end{eqnarray*}
\begin{eqnarray*} \hspace{1cm}& &\xi_i(t,\ptl_{x_1},...,\ptl_{x_n})=x_{i-1}\int_0^t(\ptl_{x_i}
+\int_0^{y_i}(\ptl_{x_{i+1}}+...+\int_0^{y_{n-2}}(\ptl_{x_{n-1}}\\
& & +y_{n-1}\ptl_{x_n}^{m_n})^{m_{n-1}}dy_{n-1}...)
^{m_{i+1}}dy_{i+1})^{m_i}dy_i\hspace{6.25cm}(5.20)\end{eqnarray*}
and
$$\xi_n(t,\ptl_{x_1},...,\ptl_{x_n})=tx_{n-1}\ptl_{x_n}^{m_n}.\eqno(5.21)$$
So
$$e^{D(t)}=e^{\xi_n(t,\ptl_{x_1},...,\ptl_{x_n})}
e^{\xi_{n-1}(t,\ptl_{x_1},...,\ptl_{x_n})}\cdots
e^{\xi_1(t,\ptl_{x_1},...,\ptl_{x_n})}.\eqno(5.22)$$

For convenience, we denote
$$k^\dg_i=\frac{k_i}{a_i},\;\;\vec
k^\dg=(k^\dg_1,...,k_n^\dg)\qquad\for\;\;\vec
k=(k_1,...,k_n)\in\mbb{N}^{\:n}.\eqno(5.23)$$ Set
$$e^{2\pi (\vec k^\dg\cdot\vec x)\sqrt{-1}}=e^{\sum_{r=1}^n2\pi
k^\dg_rx_r\sqrt{-1}}.\eqno(5.24)$$ Observe that
$$e^{D(t)}(e^{2\pi \vec (k^\dg\cdot\vec x)\sqrt{-1}})
=(\prod_{i=1}^ne^{\xi_i(t, 2\pi k^\dg_1\sqrt{-1},...,2\pi
k^\dg_n\sqrt{-1})})e^{2\pi \vec( k^\dg\cdot\vec
x)\sqrt{-1}}\eqno(5.25)$$ is a solution of (5.2) for any $\vec
k=(k_1,...,k_n)\in\mbb{N}^{\:n}$. Define
\begin{eqnarray*} \hspace{2cm}& &\phi_{\vec
k}(t,x_1,...,x_n)=\frac{1}{2}[(\prod_{i=1}^ne^{\xi_i(t, 2\pi
k^\dg_1\sqrt{-1},...,2\pi k^\dg_n\sqrt{-1})})e^{2\pi \vec
(k^\dg\cdot\vec x)\sqrt{-1}}\\ & &+(\prod_{i=1}^ne^{\xi_i(t, -2\pi
k^\dg_1\sqrt{-1},...,-2\pi k^\dg_n\sqrt{-1})})e^{-2\pi \vec
(k^\dg\cdot\vec x)\sqrt{-1}}]\hspace{4.3cm}(5.26)\end{eqnarray*}
and
\begin{eqnarray*} \hspace{2cm}& &\psi_{\vec
k}(t,x_1,...,x_n)=\frac{1}{2\sqrt{-1}}[(\prod_{i=1}^ne^{\xi_i(t,
2\pi k^\dg_1\sqrt{-1},...,2\pi k^\dg_n\sqrt{-1})})e^{2\pi \vec
(k^\dg\cdot\vec x)\sqrt{-1}}\\ & &-(\prod_{i=1}^ne^{\xi_i(t, -2\pi
k^\dg_1\sqrt{-1},...,-2\pi k^\dg_n\sqrt{-1})})e^{-2\pi \vec
(k^\dg\cdot\vec x)\sqrt{-1}}].\hspace{4.2cm}(5.27)\end{eqnarray*}
Then
$$\phi_{\vec k}(0,x_1,...,x_n)=\cos 2\pi (\vec
k^\dg\cdot\vec x),\;\;\psi_{\vec k}(0,x_1,...,x_n)=\sin 2\pi (\vec
k^\dg\cdot\vec x).\eqno(5.28)$$ By Fourier transformation theory,
we get:\psp

{\bf Theorem 5.1}. {\it The solution of the equation (5.2) subject
to (5.3) is
$$u=\sum_{\vec k\in\mbb{N}^{\:n}}(b_{\vec k}\phi_{\vec
k}(t,x_1,...,x_n)+c_{\vec k}\psi_{\vec
k}(t,x_1,...,x_n))\eqno(5.29)$$ with
$$b_{\vec k}=\frac{1}{a_1a_2\cdots a_n}\int_{-a_1}^{a_1}\cdots
\int_{-a_n}^{a_n}f(x_1,...,x_n)\cos 2\pi (\vec k^\dg\cdot\vec
x)\:dx_n\cdots dx_1\eqno(5.30)$$ and}
$$c_{\vec k}=\frac{1}{a_1a_2\cdots a_n}\int_{-a_1}^{a_1}\cdots
\int_{-a_n}^{a_n}f(x_1,...,x_n)\sin 2\pi (\vec k^\dg\cdot\vec
x)\:dx_n\cdots dx_1.\eqno(5.31)$$ \vspace{0.1cm}

{\bf Example 5.1}. Consider the case $n=2$ and $m_1=m_2=2$. So the
problem becomes
$$u_t=u_{x_1x_1}+x_1u_{x_2x_2}\eqno(5.32)$$
subject to
$$u(0,x_1,x_2)=f(x_1,x_2)\qquad\for\;\;x_i\in[-a_i,a_i].\eqno(5.33)$$
In this case,
\begin{eqnarray*}\hspace{2cm}\xi_1(t,\ptl_{x_1},\ptl_{x_2})&=&
\int_0^t(\ptl_{x_1}+y_1\ptl^2_{x_2})^2dy_1\\
&=&\int_0^t(\ptl_{x_1}^2+2y_1\ptl_{x_1}\ptl^2_{x_2}+y_1^2\ptl_{x_2}^4)dy_1\\
&=&t\ptl_{x_1}^2+t^2\ptl_{x_1}\ptl^2_{x_2}+\frac{t^3\ptl_{x_2}^4}{3}
\hspace{5.5cm}(5.34)\end{eqnarray*} and
$\xi_2(t,\ptl_{x_1},\ptl_{x_2})=tx_1\ptl_{x_2}^2$. Thus
\begin{eqnarray*}&
&e^{tx_1\ptl_{x_2}^2}e^{t\ptl_{x_1}^2+t^2\ptl_{x_1}\ptl^2_{x_2}+t^3\ptl_{x_2}^4/3}(e^{2\pi(k_1^\dg
x_1+k_2^\dg x_2)\sqrt{-1}})\\ &=&
e^{4\pi^2t(4\pi^2k_2^4t^2/3a_2^4-k_2^2x_1/a_2^2-k_1^2)}e^{2\pi(k_1x_1/a_1+k_2x_2/a_2-4\pi^2k_1k_2^2t^2/a_1a_2^2)
\sqrt{-1}}.\hspace{3.2cm}(5.35)\end{eqnarray*} Hence
\begin{eqnarray*}\hspace{1cm}\phi_{\vec
k}(t,x_1,x_2)&=&e^{4\pi^2t(4\pi^2k_2^4t^2/3a_2^4-k_2^2x_1/a_2^2-k_1^2)}\\&
&\times \cos 2\pi \left
(\frac{k_1x_1}{a_1}+\frac{k_2x_2}{a_2}-\frac{4\pi^2k_1k_2^2t^2}{a_1a_2^2}\right)\hspace{4.1cm}
(5.36)\end{eqnarray*} and
\begin{eqnarray*}\hspace{1cm}\psi_{\vec
k}(t,x_1,x_2)&=&e^{4\pi^2t(4\pi^2k_2^4t^2/3a_2^4-k_2^2x_1/a_2^2-k_1^2)}\\&
&\times \sin 2\pi \left
(\frac{k_1x_1}{a_1}+\frac{k_2x_2}{a_2}-\frac{4\pi^2k_1k_2^2t^2}{a_1a_2^2}\right).\hspace{4cm}
(5.37)\end{eqnarray*} \vspace{0.1cm}

Finally we want to study the general case. In order to find a
method of solving the problem, we exam (5.19)-(5.21) and set
$$\td\xi_n(t)=t\ptl_{x_n}^{m_n}\eqno(5.38)$$
and
\begin{eqnarray*} \hspace{1cm}& &\td\xi_i(t)=\int_0^t(\ptl_{x_i}
+\int_0^{y_i}(\ptl_{x_{i+1}}+...+\int_0^{y_{n-2}}(\ptl_{x_{n-1}}\\
& & +y_{n-1}\ptl_{x_n}^{m_n})^{m_{n-1}}dy_{n-1}...)
^{m_{i+1}}dy_{i+1})^{m_i}dy_i\hspace{6.25cm}(5.39)\end{eqnarray*}
for $i\in\ol{1,n-1}$. Then
$$\td\xi_i=\int_0^t(\ptl_{x_i}+\td\xi_{i+1}(y_i))^{m_i}dy_i\qquad\for\;\;i\in\ol{1,n-1}
\eqno(5.40)$$ and
$$\xi_1=\td\xi_1(t),\;\;\xi_i=x_{i-1}\td\xi_i(t)\qquad\for\;\;i\in\ol{2,n}.\eqno(5.41)$$

Let ${\cal T}^d$ be any tree with $n$ nodes and take
$(m_1,m_2,...,m_n)\in\mbb{N}^{\:n}$. Consider the equation:
$$u_t=(\ptl_{x_1}^{m_1}+\sum_{(\iota_i,\iota_j)\in{\cal
E}}x_i\ptl_{x_j}^{m_j})(u)\eqno(5.42)$$ subject to the condition
(5.3). Note that by the automorphism $\ast$ determined by (5.6),
$$\mbox{the Lie subalgebra of}\;\mbb{A}\;\mbox{generated
by}\;\{\ptl_{x_1}^{m_1},x_i\ptl_{x_j}^{m_j}\mid
(\iota_i,\iota_j)\in{\cal E}\}\eqno(5.43)$$ is isomorphic to the
Lie algebra ${\cal L}_0(\ol{\cal T}^d)$ for some tree diagram
$\ol{\cal T}^d$. Directly solving (5.42) would use the structure
of some abelian ideals of the Lie algebra (5.43). However, we will
use the fact of (5.40) and (5.41) and redefine $\td\xi_i(t)$ and
$\xi_i(t,\ptl_{\xi_1},...,\ptl_{x_n})$.

Recall that $\Psi$ is the set of all tips in ${\cal T}$. We set
$$\td\xi_r(t)=t\ptl_{x_r}^{m_r}\qquad\for\;\;\iota_r\in\Psi.\eqno(5.44)$$
Suppose that we have defined $\{\td\xi_s(t)\mid \iota_s\in{\cal
D}_i\}$. Set
$$\Theta_i=\{\iota_s\in{\cal N}\mid (\iota_i,\iota_s)\in{\cal
E}\}\subset{\cal D}_i.\eqno(5.45)$$ Now we define
$$\td\xi_i(t)=\int_0^t(\ptl_{x_i}+\sum_{\iota_s\in\Theta_i}\td\xi_s(y_i))^{m_i}dy_i.\eqno(5.46)$$
By induction, we have defined all
$\{\td\xi_1(t),...,\td\xi_n(t)\}$. Moreover, we let
$$\xi_1(t,\ptl_{\xi_1},...,\ptl_{x_n})=\td\xi_1(t),\;\;\xi_i(t,\ptl_{\xi_1},...,\ptl_{x_n})
=x_{p(i)}\td\xi_i(t) \eqno(5.47)$$ for $i\in\ol{2,n}$, where
$\iota_{p(i)}$ is the parent node of $\iota_i$.

Set
$$D(t)=\ptl_{x_1}^{m_1}+\sum_{(\iota_i,\iota_j)\in{\cal
E}}x_i\ptl_{x_j}^{m_j}.\eqno(5.48)$$ It can be proved that
$$e^{D(t)}=e^{\xi_n(t,\ptl_{\xi_1},...,\ptl_{x_n})}e^{\xi_{n-1}(t,\ptl_{\xi_1},...,\ptl_{x_n})}
\cdots e^{\xi_1(t,\ptl_{\xi_1},...,\ptl_{x_n})}.\eqno(5.49)$$ Thus
we have:\psp

{\bf Theorem 5.2}. {\it The solution of the equation (5.42)
subject to (5.3) is given by (5.29)-(5.31), where we replace
$\{\xi_1(t,\ptl_{\xi_1},...,\ptl_{x_n}),...,\xi_n(t,\ptl_{\xi_1},...,\ptl_{x_n})\}$
in (5.26) and (5.27) by those defined in (5.46) and (5.47).}

\vspace{1cm}

\noindent{\Large \bf References}

\hspace{0.5cm}

\begin{description}

\item[{[AKOP]}] G. Andrews, C. Krattenthaler, L. Orsina and P.
Papi, Ad-Nilpotent B-ideals in $sl(n)$ having a fixed class of
nilpotence: combinatorics and enumeration, {\it Trans. Amer. Math.
Soc.} {\bf 354} (2002), 3835-3853.

\item[{[BG1]}] J. Barros-Neto and I. M. Gel'fand, Fundamental
solutions for the Tricomi operator, {\it Duke Math. J.} {\bf 98}
(1999), 465-483.

\item[{[BG2]}] J. Barros-Neto and I. M. Gel'fand, Fundamental
solutions for the Tricomi operator II, {\it Duke Math. J.} {\bf
111} (2002), 561-584.

\item[{[CP1]}] P. Cellini and P. Papi, Ad-Nilpotent ideals of a Borel
subalgebra, {\it J. Algebra} {\bf 225} (2000), 130-141.

\item[{[CP2]}] P. Cellini and P. Papi, Enumeration of ad-nilpotent ideals of a Borel subalgebra
in type A by a class of nilpotence, {\it C. R. Acad. Math. Paris,
S\'{e}r. I, Math.} {\bf 330} (2000), 651-655.

\item[{[CP3]}] P. Cellini and P. Papi, Ad-Nilpotent ideals of a Borel
subalgebra, {\it J. Algebra} {\bf 258} (2002), 112-121.

\item[{[CP4]}] P. Cellini and P. Papi, Abelian ideals of Borel
subalgebras and affine Weyl groups, {\it Adv. Math.} {\bf 187}
(2004), 320-361.

\item[{[FS]}] G. Favre and L. Santharoubane, Nilpotent Lie algebras of
classical simple type, {\it J. Algebra} {\bf 202} (1998), 589-610.

\item[{[GK]}] M. Goze and Y. Khakimjanov, {\it Nilpotent Lie
Algebras}, ``Mathematics and its Applications'' {\bf 361}, Kluwer
Academic Publishers, Dordrecht/Boston/London. 1996.

\item[{[J1]}] N. Jacobson, Schur's theorems on commutative
matrices, {\it Bull. Amer. Math. Soc.} {\bf 50} (1944), 431-436.

\item[{[J2]}] N. Jacobson, {\it Lie algebras}, Wiley Interscience,
New York-London, 1962.

\item[{[Kv]}] V. Kac, {\it Infinite Dimensional Lie Algebras}, 3rd
Edition, Cambridge University Press, 1990.

\item[{[K1]}] B. Kostant, Eigenvalues of a Laplacian and
commutative Lie subalgebras, {\it Topology} {\bf 3} (1965),
147-159.

\item[{[K2]}] B. Kostant, The set of abelian ideals of a Borel
subalgebras, Cartan decompositions, and discrete series
representations, {\it Int. Math. Res. Not.} {\bf 5} (1998),
225-252.

\item[{[K3]}] B. Kostant, Powers of the Euler product and
commutative subalgebras of a complex Lie algebra, {\it Invent.
Math.} {\bf 158} (2004), 181-226.

\item[{[KOP]}]C. Krattenthaler, L. Orsina and P.
Papi,  Enumeration of Ad-Nilpotent B-ideals for simple Lie
algebras, {\it Adv. Appl. Math.} {\bf 28} (2002), 478-522.

\item[{[M]}] A. Malcev, Commutative subalgebras of semi-simple Lie
algebras, {\it Bull. Acad. Sci. URSS S\'{e}r. Math.} {\bf 9}
(1945), 291-300.

\item[{[MZJ]}] D. Meng, L. Zhu and C. Jiang, {\it Complete Lie
algebras}, Chinese Science Press, 2001.

\item[{[P]}] D. Panyushev, Abelian ideals of a Borel subalgebra
and long positive roots, {\it Int. Math. Res. Not.} (2003),
1889-1913.

\item[{[PR]}] D. Panyushev and G. R\"{o}hrle, Spherical orbits and
abelian ideals, {\it Adv. Math.} {\bf 159} (2001), 229-246.

\item[{[Sl1]}] L. Santharoubane, Kac-Moody Lie algebras and the
classification of nilpotent Lie algebras of maximal rank, {\it
Canad. J. Math.} {\bf 34} (1982), 1215-1239.

\item[{[Sl2]}] L. Santharoubane, Kac-Moody Lie algebras and the
universal element for the category of nilpotent Lie algebras, {\it
Math. Ann.} {\bf 263} (1983), 365-370.

\item[{[Si]}] I. Schur, Zur Theorie der vertauschbaren Matrizen,
{\it J. Reine Angrew. Math.} {\bf 130} (1905), 66-76.

\item[{[Ss]}] S. Steinberg, Applications of the Lie algebraic
formulas of Baker, Campbell, Hausdorff and Zassebhaus to the
calculation of explicit solutions of partial differential
equations, {\it J. Diff. Eq.} {\bf 26} (1977), 404-434.

\item[{[Sr]}] R. Suter, Abelian ideals in a Borel subalgebra of a
complex simple Lie algebra, {\it Invent. Math.} {\bf 156} (2004),
175-221.

\end{description}

\end{document}